\documentclass[a4paper,11pt,reqno]{amsart}

\usepackage{amssymb}

\newtheorem{theorem}[equation]{Theorem}
\newtheorem{lemma}[equation]{Lemma}
\newtheorem{corollary}[equation]{Corollary}

\numberwithin{equation}{section}

\theoremstyle{definition}

\newtheorem*{example*}{Example}

\newtheorem{examples}[equation]{Examples}
\newtheorem{remark}[equation]{Remark}
\newtheorem*{remark*}{Remark}

\newcommand{\bZ}{{\mathbb Z}}

\newcommand{\frg}{{\mathfrak g}}

\newcommand{\frh}{{\mathfrak h}}

\newcommand{\frs}{{\mathfrak s}}

\newcommand{\frv}{{\mathfrak v}}

\newcommand{\frz}{{\mathfrak z}}

\newcommand{\calS}{{\mathcal S}}

\newcommand{\calT}{{\mathcal T}}
\newcommand{\calA}{{\mathcal A}}
\newcommand{\calH}{{\mathcal H}}
\newcommand{\calU}{{\mathcal U}}

\newcommand{\calK}{{\mathcal K}}

\newcommand{\subo}{_{\bar 0}}
\newcommand{\subuno}{_{\bar 1}}

\providecommand{\espan}[1]{\text{span}\left\{ #1\right\}}

 \newcommand{\frsl}{{\mathfrak{sl}}}
 \newcommand{\frsp}{{\mathfrak{sp}}}
 \newcommand{\frso}{{\mathfrak{so}}}
 
 \newcommand{\frgl}{{\mathfrak{gl}}}

 \newcommand{\frd}{{\mathfrak{d}}}

 \DeclareMathOperator{\trace}{trace}
 \DeclareMathOperator{\ad}{ad}
 \newcommand{\der}{\mathfrak{der}}
 
 \DeclareMathOperator{\End}{End}
 
 \DeclareMathOperator{\Mat}{Mat}
 \DeclareMathOperator{\Aut}{Aut}

\newenvironment{romanenumerate}
 {\begin{enumerate}
 
 }{\end{enumerate}}

\begin{document}

\title[$S_4$-symmetry on the Tits construction]{$S_4$-symmetry on the Tits
construction of exceptional Lie algebras and superalgebras}

\author[Alberto Elduque]{Alberto Elduque$^{\star}$}
 \thanks{$^{\star}$ Supported by the Spanish Ministerio de
 Educaci\'{o}n y Ciencia
 and FEDER (MTM 2004-081159-C04-02) and by the
Diputaci\'on General de Arag\'on (Grupo de Investigaci\'on de
\'Algebra)}
 \address{Departamento de Matem\'aticas, Universidad de
Zaragoza, 50009 Zaragoza, Spain}
 \email{elduque@unizar.es}

\author[Susumu Okubo]{Susumu Okubo$^{\ast}$}
 \thanks{$^{\ast}$ Supported in part by U.S.~Department of Energy Grant No.
 DE-FG02-91ER40685.}
 \address{Department of Physics and Astronomy, University of
 Rochester, Rochester, NY 14627, USA}
 \email{okubo@pas.rochester.edu}


\date{October 20, 2006}



\begin{abstract}
The classical Tits construction provides models of the exceptional
simple Lie algebras in terms of a unital composition algebra and a
degree three simple Jordan algebra. A couple of actions of the
symmetric group $S_4$ on this construction are given. By means of
these actions, the models provided by the Tits construction are
related to models of the exceptional Lie algebras obtained from two
different types of structurable algebras. Some models of exceptional
Lie superalgebras are discussed too.
\end{abstract}

\maketitle


\section*{Introduction}

In a previous paper \cite{EOJAlgebra}, the authors have studied
those Lie algebras with an action of the symmetric group of degree
$4$, denoted by $S_4$, by automorphisms. Under some conditions,
these Lie algebras are coordinatized by the structurable algebras
introduced by Allison \cite{All78}.

The purpose of this paper is to show how a structurable algebra of
an admissible triple  appears naturally when considering an action
by automorphisms of the symmetric group $S_4$ on the classical Tits
construction of the exceptional Lie algebras \cite{Tits66}. This can
be extended to the superalgebra setting. A different $S_4$ action
will be considered too, related to the structurable algebras
consisting of a tensor product of two composition algebras. This
provides connections of the Tits construction to other models of the
exceptional Lie algebras.

The paper is structured as follows. The first section will be
devoted to show how the symmetric group $S_4$ acts by automorphisms
of the split Cayley algebra. Sections 2 and 3 will review,
respectively, the classical Tits construction \cite{Tits66} of the
exceptional Lie algebras, and the structurable algebras of
admissible triples attached to separable Jordan algebras of degree
$3$. Then Section 4 will show how to extend the action of $S_4$ on
the Cayley algebra to an action by automorphisms on the Tits
Construction. The associated coordinate algebra will be shown to be
isomorphic to the structurable algebra attached to the Jordan
algebra used in the construction. The proof involves many
computations, but the isomorphism given is quite natural. Section 5
will extend the results of the previous section to the superalgebra
setting. Here some structurable superalgebras appear as coordinate
superalgebras of the exceptional Lie superalgebras $G(3)$ and
$F(4)$. Section 6 will deal with a different action of $S_4$ on the
Tits construction. This time an action of $S_4$ by automorphisms is
given on a central simple degree $3$ Jordan algebra, and this action
is extended to an action on the Lie algebras. The associated
structurable coordinate algebra is shown to be isomorphic to the
tensor product of two unital composition algebras: the one used `on
the left' in the Tits construction, and the one that coordinatizes
the Jordan algebra involved.

All these previous results suggest a characterization of those Lie
algebras endowed with an action of $S_4$ by automorphisms in such a
way that the coordinate algebra is unital, in terms of the existence
of a subalgebra isomorphic to the three dimensional orthogonal Lie
algebra $\frso_3$, with the property that, as a module for this
subalgebra, the Lie algebra is a sum of copies of the adjoint
module, of its natural $5$ dimensional irreducible module, and of
the trivial module. This characterization is proved in Section 7.

\smallskip

In what follows, all the algebras and superalgebras considered will
be defined over a ground field $k$ of characteristic $\ne 2,3$.

\section{Composition algebras}

Let $C$ be a unital composition algebra over $k$. Thus $C$ is a
finite dimensional $k$-algebra with a nondegenerate quadratic form
$n:C\rightarrow k$ such that $n(ab)=n(a)n(b)$ for any $a,b\in C$.
Then each element $a\in c$ satisfies the degree $2$ equation:
\begin{equation}\label{eq:deg2}
a^2-t(a)a+n(a)1=0
\end{equation}
where $t(a)=n(a,1)\,\bigl(=n(a+1)-n(a)-1\bigr)$ is called the trace.
The subspace of trace zero elements will be denoted by $C^0$.

Moreover, for any $a,b\in C$, the linear map $D_{a,b}:C\rightarrow
C$ given by
\[
D_{a,b}(c)=[[a,b],c]+3(a,c,b)
\]
where $[a,b]=ab-ba$ is the commutator, and $(a,c,b)=(ac)b-a(cb)$ the
associator, is the inner derivation determined by the elements $a,b$
(see \cite[Chapter III, \S 8]{Schafer}). These derivations span the
whole Lie algebra of derivations $\der C$ of $C$. Besides, they
satisfy
\begin{equation}\label{eq:Dcyclic}
D_{a,b}=-D_{b,a},\quad D_{ab,c}+D_{bc,a}+D_{ca,b}=0,
\end{equation}
for any $a,b,c\in C$.

The dimension of $C$ is restricted to $1$, $2$, $4$ (quaternion
algebras) or $8$ (Cayley algebras), and for dimensions $2$, $4$ or
$8$ there is a unique unital composition algebra with zero divisors.
These are called split. The unique split Cayley algebras has a basis
(see \cite[Chapter 2]{ZSSS}) $\{e_1,e_2,u_0,u_1,u_2,v_0,v_1,v_2\}$
with multiplication given by:
\begin{equation}\label{eq:Csplit}
\begin{split}
&e_l^2=e_l,\ l=1,2,\ e_1e_2=0=e_2e_1,\\
&e_1u_i=u_i=u_ie_2,\ e_2v_i=v_i=v_ie_1,\ i=0,1,2,\\
&e_2u_i=0=u_ie_1,\ e_1v_i=0=v_ie_2,\ i=0,1,2,\\
&u_iu_{i+1}=v_{i+2}=-u_{i+1}u_i,\ v_iv_{i+1}=u_{i+2}=-v_{i+1}v_i,\
\textrm{indices modulo $3$,}\\
&u_i^2=0=v_i^2,\ i=0,1,2,\\
&u_iv_j=-\delta_{ij}e_1,\ v_iu_j=-\delta_{ij}e_2,\ i,j=0,1,2.
\end{split}
\end{equation}
It follows that $n(e_l)=n(u_i)=n(v_i)=0$, $l=1,2$, $i=0,1,2$, while
$n(e_1,e_2)=1=n(u_i,v_j)=\delta_{ij}$, $i,j=0,1,2$, and the unity
element is $1=e_1+e_2$.

The unique split quaternion algebra is, up to isomorphism, the
subalgebra spanned by $\{e_1,e_2,u_1,v_1\}$, which in turn is
isomorphic to the associative algebra of order $2$ matrices over
$k$: $\Mat_2(k)$. The unique split composition algebra of dimension
$2$ is the subalgebra $ke_1+ke_2$, which is isomorphic to $k\times
k$.

\medskip

The symmetric group of degree $4$, denoted by $S_4$, is generated by
the permutations
\begin{equation}\label{eq:phitaus}
\tau_1=(12)(34),\ \tau_2=(23)(14),\ \varphi=(123),\ \tau=(12),
\end{equation}
which satisfy the relations:
\[
\begin{aligned}
\tau_1\tau_2&=\tau_2\tau_1,&\varphi\tau_1&=\tau_2\varphi,&
\varphi\tau_2&=\tau_1\tau_2\varphi,\\
\tau_1\tau&=\tau\tau_1,&\tau_2\tau&=\tau\tau_2\tau_1,&
\tau\varphi&=\varphi^2\tau.
\end{aligned}
\]
The subgroup generated by $\tau_1$ and $\tau_2$ is Klein's $4$-group
$V$, the one generated by $\tau_1$, $\tau_2$ and $\varphi$ is the
alternating group $A_4$.

\medskip

Let $C$ be the split Cayley algebra and take a basis as in
\eqref{eq:Csplit}. The symmetric group $S_4$ embeds in the
automorphism group of $C$ as follows (the automorphisms of $C$ will
be denoted by the same Greek letters):
\begin{equation}\label{eq:S4actiononC}
\left\{\begin{aligned} &\textrm{$e_1$ and $e_2$ are fixed by any
element of $S_4$,}\\
&\tau_1(u_0)=u_0,\ \tau_2(u_0)=-u_0,\ \tau_1(v_0)=v_0,\
\tau_2(v_0)=-v_0,\\
&\tau_1(u_1)=-u_1,\ \tau_2(u_1)=u_1,\ \tau_1(v_1)=-v_1,\
\tau_2(v_1)=v_1,\\
&\tau_1(u_2)=-u_2,\ \tau_2(u_2)=-u_2,\ \tau_1(v_2)=-v_2,\
\tau_2(v_2)=-v_2,\\
&\varphi(u_i)=u_{i+1},\ \varphi(v_i)=v_{i+1},\ \textrm{indices
modulo
$3$,}\\
&\tau(u_0)=-u_0,\ \tau(u_1)=-u_2,\ \tau(u_2)=-u_1,\\
&\tau(v_0)=-v_0,\ \tau(v_1)=-v_2,\ \tau(v_2)=-v_1.
\end{aligned}\right.
\end{equation}

The action of Klein's $4$ group $V$ gives a grading of $C$ over
$\bZ_2\times\bZ_2$:
\begin{equation}\label{eq:CZ2Z2grading}
C=C_{(\bar 0,\bar 0)}\oplus C_{(\bar 1,\bar 0)}\oplus C_{(\bar
0,\bar 1)}\oplus C_{(\bar 1,\bar 1)},
\end{equation}
where
\begin{equation}\label{eq:CZ2Z2gradingbis}
\begin{aligned}
C_{(\bar 0,\bar 0)}&=ke_1+ke_2,&C_{(\bar 1,\bar 0)}&=ku_0+kv_0,\\
C_{(\bar 0,\bar 1)}&=ku_1+kv_1,&C_{(\bar 1,\bar 1)}&=ku_2+kv_2.
\end{aligned}
\end{equation}

Any automorphism $\psi$ of $C$ induces an automorphism of $\der C$:
$d\mapsto \psi d\psi^{-1}$. Note that one has $D_{1,a}=0$ for any
$a\in C$, so that
\[
\der
C=D_{C,C}=D_{C^0,C^0}=D_{U,V}+D_{U,U}+D_{V,V}+D_{e_1-e_2,U}+D_{e_2-e_1,V},
\]
where $U$ (respectively $V$) denotes here the span of the $u_i$'s
(resp. $v_i$'s). But, because of \eqref{eq:Dcyclic}
\[
D_{U,U}=D_{V^2,U}\subseteq D_{V,UV+VU}\subseteq
D_{ke_1+ke_2,V}=D_{e_2-e_1,V},
\]
and, similarly, $D_{V,V}\subseteq D_{e_1-e_2,U}$. Since the
decomposition $C=(ke_1+ke_2)\oplus U\oplus V$ is a grading of $C$
over $\bZ_3$, which induces a grading of $\der C$, it follows that
\[
\der C=D_{U,V}\oplus D_{e_1-e_2,U}\oplus D_{e_2-e_1,V}
\]
is the associated grading of $\der C$ over $\bZ_3$.

The action of Klein's $4$ group on $\der C$ produces an associated
grading on $\der C$:
\begin{equation}\label{eq:derCZ2Z2grading}
\der C=(\der C)_{(\bar 0,\bar 0)}\oplus (\der C)_{(\bar 1,\bar
0)}\oplus (\der C)_{(\bar 0,\bar 1)}\oplus (\der C)_{(\bar 1,\bar
1)},
\end{equation}
Let us compute $(\der C)_{(\bar 1,\bar 0)}=D_{C_{(\bar 0,\bar
0)},C_{(\bar 1,\bar 0)}}+D_{C_{(\bar 0,\bar 1)},C_{(\bar 1,\bar
1)}}$. Because of \eqref{eq:Dcyclic}
\begin{equation}\label{eq:Du1u2}
\begin{split}
D_{u_1,u_2}&=D_{v_2v_0,u_2}=-D_{v_0u_2,v_2}-D_{u_2v_2,v_0}\\
&=D_{e_1,v_0}=-D_{\frac{1}{2}-e_1,v_0}=-\frac{1}{2}D_{e_2-e_1,v_0},
\end{split}
\end{equation}
and, with the same argument,
$D_{v_1,v_2}=-\frac{1}{2}D_{e_1-e_2,u_0}$. Therefore,
\begin{equation}\label{eq:derC10}
(\der C)_{(\bar 1,\bar 0)}=\espan{D_{e_1-e_2,u_0},D_{e_2-e_1,v_0},
D_{u_1,v_2},D_{v_1,u_2}}.
\end{equation}
Using the multiplication table in \eqref{eq:Csplit}, the next
equations follow:
\[
\begin{split}
D_{u_i,v_j}(e_l)&=[[u_i,v_j],e_l]+3(u_i,e_l,v_j)=0,\\
D_{e_1-e_2,u_i}(e_1)&=[[e_1-e_2,u_i],e_1]+3(e_1-e_2,e_1,u_i)=2[u_i,e_1]=-2u_i,\\
D_{u_i,v_{i+1}}(u_i)&=3(u_i,u_i,v_{i+1})=0,\\
D_{u_i,v_{i+1}}(v_i)&=3(u_i,v_i,v_{i+1})=3v_{i+1},
\end{split}
\]
for any $i,j=0,1,2$ (indices modulo $3$), and by symmetry
$e_1\leftrightarrow e_2$, $u_i\leftrightarrow v_i$, one has also
\[
D_{v_i,u_{i+1}}(v_i)=0,\quad D_{v_i,u_{i+1}}(u_i)=3u_{i+1}.
\]
From here it follows that the elements in \eqref{eq:derC10} are
linearly independent, so that $\dim (\der C)_{(\bar 1,\bar 0)}=4$.

\section{Tits construction}

Some results in \cite[Sections 3 and 4]{BZ96} (see also
\cite{Tits66} and \cite{BE03}) will be reviewed in this section.

Let $C$ be a unital composition algebra over the ground field $k$
with norm $n$ and trace $t$. Let $J$ be a unital Jordan algebra with
a \emph{normalized trace} $t_J:J\rightarrow k$. That is, $t_J$ is a
linear map such that $t_J(1)=1$ and
$t_J\bigl((xy)z\bigr)=t_J\bigl(x(yz)\bigr)$ for any $x,y,z\in J$.
Then $J=k1\oplus J^0$, where $J^0=\{x\in J: t_J(x)=0\}$. For $x,y\in
J^0$,
\begin{equation}\label{eq:Jxy}
xy=t_J(xy)1+x*y,
\end{equation}
where $x*y=xy-t_J(xy)1$ gives a commutative multiplication on $J^0$.
For $x,y\in J$, the linear map $d_{x,y}:J\rightarrow J$ defined by
\begin{equation}\label{eq:dxy}
d_{x,y}(z)=x(yz)-y(xz),
\end{equation}
is the inner derivation of $J$ determined by the elements $x$ and
$y$. Since $d_{1,x}=0$ for any $x$, it is enough to deal with the
inner derivations $d_{x,y}$, with $x,y\in J^0$.

Given $C$ and $J$ as before, consider the space
\begin{equation}\label{eq:TCJ}
\calT(C,J)=\der C\oplus \bigl(C^0\otimes J^0\bigr)\oplus d_{J,J}
\end{equation}
(unadorned tensor products are always considered over $k$), with the
anticommutative multiplication $[.,.]$ specified by:
\begin{equation}\label{eq:TCJproduct}
\begin{split}
\bullet&\ \textrm{$\der C$ and $d_{J,J}$ are Lie subalgebras,}\\
\bullet&\ [\der C,d_{J,J}]=0,\\
\bullet&\ [D,a\otimes x]=D(a)\otimes x,\ [d,a\otimes x]=a\otimes
d(x),\\
\bullet&\ [a\otimes x,b\otimes y]=t_J(xy)D_{a,b}+\bigl([a,b]\otimes
x*y\bigr)+2t(ab)d_{x,y},
\end{split}
\end{equation}
for all $D\in \der C$, $d\in d_{J,J}$, $a,b\in C^0$, and $x,y\in
J^0$. Here the bracket $[.,.]$ follows the conventions in
\cite[(1.4)]{BE03}.

The conditions for $\calT(C,J)$ to be a Lie algebra are the
following:
\begin{equation}\label{eq:TCJLie}
\begin{split}
\textrm{(i)}&\  \displaystyle{\sum_{\circlearrowleft}
 t\bigl([a_{1}, a_{2}] a_{3}\bigr)\,
d_{(x_1 * x_2), x_3}}=0,\\[6pt]
\textrm{(ii)}&\  \displaystyle{\sum_{\circlearrowleft}
 t_J\bigl( (x_1 * x_2) x_{3}\bigr)
\,D_{[a_1, a_2], a_3}}=0,\\[6pt]
\textrm{(iii)}&\ \displaystyle{\sum_{\circlearrowleft}
 \Bigl(D_{a_1,a_2}(a_3) \otimes t_J\bigl(x_1
x_2\bigr) x_3} + [[a_1, a_2],a_3] \otimes (x_1 * x_2)* x_3\\[-6pt]
  &\qquad\qquad\qquad\qquad +2
t(a_1 a_2) a_3\otimes d_{x_1, x_2}(x_3)\Bigr)=0
\end{split}
\end{equation}
for any $a_1,a_2,a_{3} \in C^{0}$ and any $x_1,x_2,x_3 \in J^0$. The
notation ``$\displaystyle{\sum_\circlearrowleft}$'' indicates
summation over the cyclic permutation of the variables.

These conditions appear in \cite[Proposition 1.5]{BE03}, but there
they are stated in the more general setting of superalgebras, a
setting we will deal with later on. In particular, these conditions
are fulfilled if $J$ is a separable Jordan algebra of degree three
over $k$ and $t_J=\frac{1}{3}T$, where $T$ denotes the generic trace
of $J$ (see for instance \cite{JacobsonJordan}).

\section{The algebra $\bigl(\calA(J),-\bigr)$}\label{se:AJ}

Let $J$ be a unital Jordan algebra over $k$ with a normalized trace
$t_J$ as in the previous section. For any $x,y\in J$, consider the
new commutative product on $J$ defined by
\begin{equation}\label{eq:cross}
x\times y=2xy-3t_J(x)y-3t_J(y)x+\bigl(9t_J(x)t_J(y)-3t_J(xy)\bigr)1,
\end{equation}
for any $x,y\in J$. Note that for any $x,y\in J^0$ the following
holds:
\begin{subequations}\label{eq:severalcross}
\begin{align}
1\times 1&= 2,\label{eq:1x1}\\
1\times x&= -x,\label{eq:1xx}\\
x\times y&= 2xy-3t_J(xy)1=2x*y-t_J(xy)1.\label{eq:xxy}
\end{align}
\end{subequations}

In case $J$ is a separable Jordan algebra of degree $3$ with generic
trace $T$ and generic norm $N$, and with $t_J=\frac{1}{3}T$, this is
the cross product considered in \cite[page 148]{All78},
corresponding to the admissible triple $(T,N,N)$ on the pair
$(J,J)$. Now, as in \cite{All78}, consider the space
\[
\calA(J)=\left\{\begin{pmatrix} \alpha &x\\ y&\beta\end{pmatrix}:
\alpha,\beta\in k,\ x,y\in J\right\},
\]
with multiplication given by
\[
\begin{pmatrix} \alpha &x\\ y&\beta\end{pmatrix}
 \begin{pmatrix} \alpha' &x'\\ y'&\beta'\end{pmatrix} =
 \begin{pmatrix} \alpha\alpha'+3t_J(xy') &
                  \alpha x'+\beta'x+y\times y'\\
                  \alpha' y+\beta y'+x\times x' &
                  \beta\beta'+3t_J(yx')\end{pmatrix}
\]
for any $\alpha,\beta,\alpha',\beta'\in k$, and $x,y,x',y'\in J$. It
follows that the map:
\[
\begin{pmatrix} \alpha &x\\ y&\beta\end{pmatrix}\xrightarrow{-}
\begin{pmatrix} \beta &x\\ y&\alpha\end{pmatrix}
\]
is an involution (involutive antiautomorphism) of $\calA(J)$.
Besides, $\calA(J)$ is unital with $1=\left(\begin{smallmatrix}
1&0\\ 0&1\end{smallmatrix}\right)$.

\smallskip

\begin{remark}\label{re:cubicadmissible}
Let $A$ be a commutative algebra endowed with a cubic form
$N:A\rightarrow k$ such that $(x^2)^2=N(x)x$ for any $x\in A$. These
algebras have been called \emph{admissible cubic algebras} in
\cite{EOMathZ}, where the relationships of these algebras to Jordan
algebras have been studied. Then there is a symmetric associative
bilinear form $\langle .\vert.\rangle$ on $A$, called the
\emph{trace form}, such that $N(x)=\langle x\vert x^2\rangle$ for
any $x\in A$. If $N\ne 0$, then $\langle .\vert .\rangle$ is
uniquely determined. Then the trilinear form given by
$N(x,y,z)=6\langle x\vert yz\rangle$ satisfies $N(x,x,x)=6N(x)$ for
any $x$, and hence $\bigl(3\langle .\vert .\rangle,N,N\bigr)$ is an
admissible triple on $(A,A)$ in the sense of \cite[page 148]{All78},
with associated cross product $x\times y=2xy$ for any $x,y\in A$, so
$x^\sharp =\frac{1}{2}x\times x=x^2$. Then, as before, the linear
space
\[
\calA(A)=\left\{\begin{pmatrix} \alpha &x\\ y&\beta\end{pmatrix}:
\alpha,\beta\in k,\ x,y\in A\right\},
\]
with multiplication given by
\begin{equation}\label{eq:TNNx}
\begin{pmatrix} \alpha &x\\ y&\beta\end{pmatrix}
 \begin{pmatrix} \alpha' &x'\\ y'&\beta'\end{pmatrix} =
 \begin{pmatrix} \alpha\alpha'+3\langle x\vert y'\rangle &
                  \alpha x'+\beta'x+2yy'\\
                  \alpha' y+\beta y'+2xx' &
                  \beta\beta'+3\langle y\vert x'\rangle\end{pmatrix}
\end{equation}
for any $\alpha,\beta,\alpha',\beta'\in k$, and $x,y,x',y'\in A$, is
a structurable algebra, whose involution is given by
\[
\begin{pmatrix} \alpha &x\\ y&\beta\end{pmatrix}\xrightarrow{-}
\begin{pmatrix} \beta &x\\ y&\alpha\end{pmatrix}.
\]
\end{remark}

\section{An $S_4$-action on the Tits construction}

Let $C$ be the split Cayley algebra over $k$ and let $J$ be any
unital Jordan algebra with a normalized trace $t_J$ such that
$\calT=\calT(C,J)$ is a Lie algebra. Then the action of the
symmetric group $S_4$ as automorphisms of $C$ in
\eqref{eq:S4actiononC} extends to an action by automorphisms on
$\calT=\der C\oplus\bigl(C^0\otimes J^0\bigr)\oplus d_{J,J}$ given
by:
\begin{equation}\label{eq:S4actiononTCJ}
\psi\bigl(D+(a\otimes x)+d\bigr)=\psi D\psi^{-1}+(\psi(a)\otimes
x)+d,
\end{equation}
for any $\psi\in S_4$, $D\in \der C$, $d\in d_{J,J}$, $a\in C^0$ and
$x\in J^0$.

As in \eqref{eq:CZ2Z2grading}, the action of Klein's $4$-group
induces a grading over $\bZ_2\times\bZ_2$:
\begin{equation}\label{eq:TZ2Z2grading}
\calT=\calT(C,J)=\calT_{(\bar 0,\bar 0)}\oplus \calT_{(\bar 1,\bar
0)}\oplus \calT_{(\bar 0,\bar 1)}\oplus \calT_{(\bar 1,\bar 1)}.
\end{equation}
Under these circumstances (see \cite[Section 2]{EOJAlgebra}), the
subspace $\calT_{(\bar 1,\bar 0)}=\{ X\in\calT: \tau_1(X)=X,\
\tau_2(X)=-X\}$ becomes an algebra with involution by means of:
\begin{equation}\label{eq:T10product}
\begin{split}
&X\cdot Y=-\tau\bigl([\varphi(X),\varphi^2(Y)]\bigr),\\
&\bar X=-\tau(X),
\end{split}
\end{equation}
for any $X,Y\in\calT_{(\bar 1,\bar 0)}$. Here $\varphi=(123)$ and
$\tau=(12)$ as in \eqref{eq:phitaus}.

\medskip

The purpose of this section is to prove the following result:

\begin{theorem}\label{th:T10AJ}
Let $C$ be the split Cayley algebra over $k$ and let $J$ be a unital
Jordan algebra with a normalized trace $t_J$ over $k$ such that the
Tits algebra $\calT=\calT(C,J)$ is a Lie algebra. Then the algebra
with involution $\bigl(\calT_{(\bar 1,\bar 0)},\cdot,-\bigr)$ is
isomorphic to the algebra $\bigl(\calA(J),-\bigr)$ defined in
Section \ref{se:AJ}.
\end{theorem}

Before going into the proof of this result, let us note the next
consequence, which follows immediately from \cite[Theorem
2.9]{EOJAlgebra} and the fact that the algebra $\calA(J)$ is unital:

\begin{corollary}\label{co:AJstructurable}
Under the conditions of Theorem \ref{th:T10AJ}, the algebra with
involution $\bigl(\calA(J),-\bigr)$ is a structurable algebra.
\end{corollary}

\medskip

To prove Theorem \ref{th:T10AJ}, first note that, because of
\eqref{eq:S4actiononTCJ}, \eqref{eq:CZ2Z2gradingbis}, and
\eqref{eq:derC10}, one has
\[
\begin{split}
\calT_{(\bar 1,\bar 0)}&=\bigl\{ X\in\calT: \tau_1(X)=X,\
\tau_2(X)=-X\bigr\}\\
 &=(\der C)_{(\bar 1,\bar 0)}\oplus\bigl( C_{(\bar 1,\bar 0)}\otimes
 J^0\bigr)\\
 &=\espan{D_{u_1,v_2},D_{v_1,u_2},D_{e_1-e_2,u_0},D_{e_2-e_1,v_0}}
 \oplus\bigl(u_0\otimes J^0\bigr)\oplus \bigl(v_0\otimes J^0\bigr).
\end{split}
\]
Consider the bijective linear map
\[
\begin{split}
\Phi:\calT_{(\bar 1,\bar 0)}&\longrightarrow \calA(J)\\
D_{v_1,u_2}&\mapsto \begin{pmatrix} 3&0\\ 0&0\end{pmatrix}\\
D_{u_1,v_2}&\mapsto \begin{pmatrix} 0&0\\ 0&3\end{pmatrix}\\
D_{e_1-e_2,u_0}&\mapsto \begin{pmatrix} 0&2\\ 0&0\end{pmatrix}\\
D_{e_2-e_1,v_0}&\mapsto \begin{pmatrix} 0&0\\ 2&0\end{pmatrix}\\
u_0\otimes x&\mapsto \begin{pmatrix} 0&x\\ 0&0\end{pmatrix}\\
v_0\otimes x&\mapsto \begin{pmatrix} 0&0\\ x&0\end{pmatrix}
\end{split}
\]
where $x\in J^0$. To prove Theorem \ref{th:T10AJ} it is enough to
prove that $\Phi$ is a homomorphism of algebras with involution.
Note that there is the order $2$ automorphism of $C$ given by
$e_1\leftrightarrow e_2$, $u_i\leftrightarrow v_i$, which extends to
an order $2$ automorphism $\epsilon$ of $\calT(C,J)$. On the other
hand, there is the natural order $2$ automorphism of $\calA(J)$
given by $\left(\begin{smallmatrix} \alpha&x\\
y&\beta\end{smallmatrix}\right)\overset{\delta}\leftrightarrow \left(\begin{smallmatrix} \beta&y\\
x&\alpha\end{smallmatrix}\right)$ which satisfies
$\Phi\epsilon=\delta\Phi$. This simplifies the number of
computations to be done. Thus, it is enough to check that
$\Phi(X\cdot Y)=\Phi(X)\Phi(Y)$ for the following pairs $X$, $Y$ in
$\calT_{(\bar 1,\bar 0)}$:

\smallskip

\begin{romanenumerate}

\item $X=Y=D_{v_1,u_2}$: Here
\[
\begin{split}
D_{v_1,u_2}\cdot D_{v_1,u_2}
 &=-\tau\Bigl([\varphi(D_{v_1,u_2}),\varphi^2(D_{v_1,u_2})]\Bigr)\\
 &=-\tau\Bigl([D_{v_2,u_0},D_{v_0,u_1}]\Bigr)\\
 &=-\tau\Bigl(D_{D_{v_2,u_0}(v_0),u_1}+D_{v_0,D_{v_2,u_0}(u_1)}\Bigr).
\end{split}
\]
But $D_{v_2,u_0}(v_0)=3(v_2,v_0,u_0)=-3v_2$ and
$D_{v_2,u_0}(u_1)=3(v_2,u_1,u_0)=0$, so that
\[
\begin{split}
D_{v_1,u_2}\cdot D_{v_1,u_2}
 &=-\tau\bigl(-3D_{v_2,u_1}\bigr)=-\tau\bigl(3D_{u_1,v_2}\bigr)\\
 &=-3D_{\tau(u_1),\tau(v_2)}=-3D_{u_2,v_1}=3D_{v_1,u_2}.
\end{split}
\]
Therefore,
\[
\Phi\bigl(D_{v_1,u_2}\cdot
D_{v_1,u_2}\bigr)=3\Phi\bigl(D_{v_1,u_2}\bigr)=\left(\begin{smallmatrix}
9&0\\ 0&0\end{smallmatrix}\right),
\]
while
\[
\Phi\bigl(D_{v_1,u_2}\bigr)\Phi\bigl(D_{v_1,u_2}\bigr)=
\left(\begin{smallmatrix} 3&0\\ 0&0\end{smallmatrix}\right)
\left(\begin{smallmatrix} 3&0\\ 0&0\end{smallmatrix}\right)=
\left(\begin{smallmatrix} 9&0\\ 0&0\end{smallmatrix}\right).
\]

\smallskip
\item $X=D_{v_1,u_2}$, $Y=D_{u_1,v_2}$: This is easier since
\[
[\varphi(D_{v_1,u_2}),\varphi^2(D_{u_1,v_2})]=[D_{v_2,u_0},D_{u_0,v_1}]=0,
\]
as $D_{v_2,u_0}(u_0)=0=D_{v_2,u_0}(v_1)$. Hence both $X\cdot Y$ and
$\Phi(X)\Phi(Y)$ are $0$.

\smallskip
\item $X=D_{v_1,u_2}$, $Y=\frac{1}{2}D_{e_1-e_2,u_0}=D_{e_1,u_0}$:
Here
\[
[\varphi(D_{v_1,u_2}),\varphi^2(D_{e_1,u_0})]=[D_{v_2,u_0},D_{e_1,u_2}]=3D_{e_1,u_0},
\]
because $D_{v_2,u_0}(e_1)=0$ and
$D_{v_2,u_0}(u_2)=3(v_2,u_2,u_0)=-3v_2v_1=3u_0$. Thus,
\[
X\cdot
Y=-\tau\bigl(3D_{e_1,u_0}\bigr)=-3D_{\tau(e_1),\tau(u_0)}=3D_{e_1,u_0}=3Y,
\]
while
\[
\Phi(X)\Phi(Y)=\left(\begin{smallmatrix} 3&0\\
0&0\end{smallmatrix}\right) \left(\begin{smallmatrix} 0&1\\
0&0\end{smallmatrix}\right)= \left(\begin{smallmatrix} 0&3\\
0&0\end{smallmatrix}\right)=3\Phi(Y).
\]

\smallskip
\item $X=D_{v_1,u_2}$, $Y=\frac{1}{2}D_{e_2-e_1,v_0}=D_{e_2,v_0}$:
In this case $\Phi(X)\Phi(Y)=0$, but also
\[
[\varphi(D_{v_1,u_2}),\varphi^2(D_{e_2,v_0})]=[D_{v_2,u_0},D_{e_2,v_2}]=0,
\]
as $D_{v_2,u_0}(e_2)=0=D_{v_2,u_0}(v_2)$, so $X\cdot Y=0$ too.

\smallskip
\item $X=D_{v_1,u_2}$, $Y=u_0\otimes x$, with $x\in J^0$: Here
\[
\begin{split}
[\varphi(D_{v_1,u_2}),\varphi^2(u_0\otimes
x)]&=[D_{v_2,u_0},u_2\otimes x]=D_{v_2,u_0}(u_2)\otimes x\\
&=3(v_2,u_2,u_0)\otimes x=3u_0\otimes x.
\end{split}
\]
Hence, $X\cdot
Y=-\tau(3u_0\otimes x)=3u_0\otimes x=3Y$, while
\[
\Phi(X)\Phi(Y)=\left(\begin{smallmatrix} 3&0\\
0&0\end{smallmatrix}\right) \left(\begin{smallmatrix} 0&x\\
0&0\end{smallmatrix}\right)= \left(\begin{smallmatrix} 0&3x\\
0&0\end{smallmatrix}\right)=3\Phi(Y).
\]

\smallskip
\item $X=D_{v_1,u_2}$, $Y=v_0\otimes x$, with $x\in J^0$: Then
\[
[\varphi(D_{v_1,u_2}),\varphi^2(v_0\otimes
x)]=[D_{v_2,u_0},v_2\otimes x]=D_{v_2,u_0}(v_2)\otimes x=0,
\]
and
both $X\cdot Y$ and $\Phi(X)\Phi(Y)=0$.

\smallskip
\item $X=D_{e_1,u_0}$, $Y=D_{v_1,u_2}$:
$[\varphi(X),\varphi^2(Y)]=[D_{e_1,u_1},D_{v_0,u_1}]=0$, as
$D_{e_1,u_1}(v_0)=0=D_{e_1,u_1}(u_1)$, so $X\cdot Y=0$, but also
\[
\Phi(X)\Phi(Y)=\left(\begin{smallmatrix} 0&1\\
0&0\end{smallmatrix}\right) \left(\begin{smallmatrix} 3&0\\
0&0\end{smallmatrix}\right)= 0.
\]

\smallskip
\item $X=D_{e_1,u_0}$, $Y=D_{u_1,v_2}$: One has
\[
\begin{split}
[\varphi(X),\varphi^2(Y)]&=[D_{e_1,u_1},D_{u_0,v_1}]\\
 &=-D_{D_{u_0,v_1}(e_1),u_1}
-D_{e_1,D_{u_0,v_1}(u_1)}\\
 &=-D_{e_1,3(u_0,u_1,v_1)}=3D_{e_1,u_0},
 \end{split}
\]
so that $X\cdot Y=-\tau(3D_{e_1,u_0})=3D_{e_1,u_0}=3X$, while
\[
\Phi(X)\Phi(Y)=\left(\begin{smallmatrix} 0&1\\
0&0\end{smallmatrix}\right) \left(\begin{smallmatrix} 0&0\\
0&3\end{smallmatrix}\right)= \left(\begin{smallmatrix} 0&3\\
0&0\end{smallmatrix}\right)=3\Phi(X).
\]

\smallskip
\item $X=D_{e_1,u_0}$, $Y=D_{e_1,u_0}$: Here
\[
\begin{split}
[\varphi(X),\varphi^2(Y)]&=[D_{e_1,u_1},D_{e_1,u_2}]\\
 &=D_{D_{e_1,u_1}(e_1),u_2}+D_{e_1,D_{e_1,u_1}(u_2)}\\
 &=D_{-u_1,u_2}+D_{e_1,-v_0}=-2D_{e_1,v_0}=2D_{e_2,v_0}
\end{split}
\]
(see
\eqref{eq:Du1u2}). Hence $X\cdot
Y=-\tau\bigl(2D_{e_2,v_0}\bigr)=2D_{e_2,v_0}$ and $\Phi(X\cdot Y)=
\left(\begin{smallmatrix} 0&0\\
2&0\end{smallmatrix}\right)$, while
\[
\Phi(X)\Phi(Y)=\left(\begin{smallmatrix} 0&1\\
0&0\end{smallmatrix}\right) \left(\begin{smallmatrix} 0&1\\
0&0\end{smallmatrix}\right)= \left(\begin{smallmatrix} 0&0\\
1\times 1&0\end{smallmatrix}\right)=\left(\begin{smallmatrix} 0&0\\
2&0\end{smallmatrix}\right),
\]
because of \eqref{eq:1x1}.

\smallskip
\item $X=D_{e_1,u_0}$, $Y=D_{e_2,v_0}$: Then
\[
\begin{split}
[\varphi(X),\varphi^2(Y)]&=[D_{e_1,u_1},D_{e_2,v_2}]\\
 &=
D_{D_{e_1,u_1}(e_2),v_2}+D_{e_2,D_{e_1,u_1}(v_2)}\\
 &=
D_{u_1,v_2}+0=D_{u_1,v_2}.
\end{split}
\]
Thus, $X\cdot
Y=-\tau\bigl(D_{u_1,v_2}\bigr)=-D_{\tau(u_1),\tau(v_2)}=-D_{u_2,v_1}=D_{v_1,u_2}$,
so $\Phi(X\cdot Y)=\left(\begin{smallmatrix} 3&0\\
0&0\end{smallmatrix}\right)$, while
\[
\Phi(X)\Phi(Y)=\left(\begin{smallmatrix} 0&1\\
0&0\end{smallmatrix}\right) \left(\begin{smallmatrix} 0&0\\
1&0\end{smallmatrix}\right)= \left(\begin{smallmatrix} 3&0\\
0&0\end{smallmatrix}\right).
\]

\smallskip
\item $X=D_{e_1,u_0}$, $Y=u_0\otimes x$, $x\in J^0$: In this case
\[
[\varphi(X),\varphi^2(Y)]=[D_{e_1,u_1},u_2\otimes x]=
D_{e_1,u_1}(u_2)\otimes x=-v_0\otimes x,
\]
so $X\cdot
Y=-\tau(-v_0\otimes x)=-v_0\otimes x$ and $\Phi(X\cdot Y)
=\left(\begin{smallmatrix} 0&0\\
-x&0\end{smallmatrix}\right)$, while
\[
\Phi(X)\Phi(Y)=\left(\begin{smallmatrix} 0&1\\
0&0\end{smallmatrix}\right) \left(\begin{smallmatrix} 0&x\\
0&0\end{smallmatrix}\right)= \left(\begin{smallmatrix} 0&0\\
1\times x&0\end{smallmatrix}\right)=\left(\begin{smallmatrix} 0&0\\
-x&0\end{smallmatrix}\right),
\]
because of \eqref{eq:1xx}.

\smallskip
\item $X=D_{e_1,u_0}$, $Y=v_0\otimes x$, $x\in J^0$: For these $X$
and $Y$
\[
[\varphi(X),\varphi^2(Y)]=[D_{e_1,u_1},v_2\otimes x]=
D_{e_1,u_1}(v_2)\otimes x=0,
\]
so $X\cdot Y=0$. Here
$\Phi(x)\Phi(Y)=0$ too, since $t_J(x)=0$.

\smallskip
\item $X=u_0\otimes x$, $Y=D_{v_1,u_2}$, $x\in J^0$: Here
\[
[\varphi(X),\varphi^2(Y)]=[u_1\otimes
x,D_{v_0,u_1}]=-D_{v_0,u_1}(u_1)\otimes x=0,
\]
so $X\cdot Y=0$ and
also
\[
\Phi(X)\Phi(Y)=\left(\begin{smallmatrix} 0&x\\
0&0\end{smallmatrix}\right) \left(\begin{smallmatrix} 3&0\\
0&0\end{smallmatrix}\right)= 0.
\]

\smallskip
\item $X=u_0\otimes x$, $Y=D_{u_1,v_2}$, $x\in J^0$: Now
\[
\begin{split}
[\varphi(X),\varphi^2(Y)]&=
 [u_1\otimes x,D_{u_0,v_1}]=-D_{u_0,v_1}(u_1)\otimes x\\
 &=-3(u_0,u_1,v_1)\otimes x=3u_0\otimes x,
\end{split}
\]
so $X\cdot Y=-\tau(3u_0\otimes x)=3u_0\otimes
x=3X$, while
\[
\Phi(X)\Phi(Y)=\left(\begin{smallmatrix} 0&x\\
0&0\end{smallmatrix}\right) \left(\begin{smallmatrix} 0&0\\
0&3\end{smallmatrix}\right)= \left(\begin{smallmatrix} 0&3x\\
0&0\end{smallmatrix}\right)=3\Phi(X).
\]

\smallskip
\item $X=u_0\otimes x$, $Y=D_{e_1,u_0}$, $x\in J^0$: Here
\[
[\varphi(X),\varphi^2(Y)]
 =[u_1\otimes x,D_{e_1,u_2}]=-D_{e_1,u_2}(u_1)\otimes x=-v_0\otimes x,
\]
so $X\cdot Y=-\tau(-v_0\otimes x)=-v_0\otimes x$ and $\Phi(X\cdot
Y)=
\left(\begin{smallmatrix} 0&0\\
-x&0\end{smallmatrix}\right)$, while
\[
\Phi(X)\Phi(Y)=\left(\begin{smallmatrix} 0&x\\
0&0\end{smallmatrix}\right) \left(\begin{smallmatrix} 0&1\\
0&0\end{smallmatrix}\right)= \left(\begin{smallmatrix} 0&0\\
x\times 1&0\end{smallmatrix}\right)=
\left(\begin{smallmatrix} 0&0\\
-x&0\end{smallmatrix}\right),
\]
because of \eqref{eq:1xx}.

\smallskip
\item $X=u_0\otimes x$, $Y=D_{e_2,v_0}$, $x\in J^0$: In this case
\[
[\varphi(X),\varphi^2(Y)]=[u_1\otimes
x,D_{e_2,v_2}]=-D_{e_2,v_2}(u_1)\otimes x=0,
\]
so  $\Phi(X\cdot
Y)=0$,
 while
$
\Phi(X)\Phi(Y)=\left(\begin{smallmatrix} 0&x\\
0&0\end{smallmatrix}\right) \left(\begin{smallmatrix} 0&0\\
1&0\end{smallmatrix}\right)= 0$, as $t_J(x)=0$.

\smallskip
\item $X=u_0\otimes x$, $Y=u_0\otimes y$, $x,y\in J^0$: Here
\[
\begin{split}
[\varphi(X),\varphi^2(Y)]
 &=[u_1\otimes x,u_2\otimes y] \\
 &=t_J(xy)D_{u_1,u_2}+[u_1,u_2]\otimes x*y+2t(u_1u_2)d_{x,y}\\
 &=-t_J(xy)D_{e_2,v_0}+2v_0\otimes x*y
\end{split}
\]
(see \eqref{eq:Du1u2}). Hence
\[
X\cdot Y=-\tau\bigl(-t_J(xy)D_{e_2,v_0}+2v_0\otimes
x*y\bigr)=-t_J(xy)D_{e_2,v_0}+2v_0\otimes x*y,
\]
and
\[\Phi(X\cdot Y)=
\left(\begin{smallmatrix} 0&0\\
-t_J(xy)+2x*y&0\end{smallmatrix}\right)=\left(\begin{smallmatrix} 0&0\\
x\times y&0\end{smallmatrix}\right),
\]
using \eqref{eq:xxy}, while
\[
\Phi(X)\Phi(Y)=\left(\begin{smallmatrix} 0&x\\
0&0\end{smallmatrix}\right) \left(\begin{smallmatrix} 0&y\\
0&0\end{smallmatrix}\right)= \left(\begin{smallmatrix} 0&0\\
x\times y&0\end{smallmatrix}\right).
\]

\smallskip
\item $X=u_0\otimes x$, $Y=v_0\otimes y$, $x,y\in J^0$: In this
final case
\[
[\varphi(X),\varphi^2(Y)]=[u_1\otimes x,v_2\otimes
y]=t_J(xy)D_{u_1,v_2}+[u_1,v_2]\otimes x*y=t_J(xy)D_{u_1,v_2},
\]
so
\[
\begin{split}
X\cdot Y
 &=-\tau\bigl(t_J(xy)D_{u_1,v_2}\bigr)=-t_J(xy)D_{\tau(u_1),\tau(v_2)}\\
 &=-t_J(xy)D_{u_2,v_1}=t_J(xy)D_{v_1,u_2},
\end{split}
\]
and $\Phi(X\cdot Y)=\left(\begin{smallmatrix} 3t_J(xy)&0\\
0&0\end{smallmatrix}\right)$, while
\[
\Phi(X)\Phi(Y)=\left(\begin{smallmatrix} 0&x\\
0&0\end{smallmatrix}\right) \left(\begin{smallmatrix} 0&0\\
y&0\end{smallmatrix}\right)= \left(\begin{smallmatrix} 3t_J(xy)&0\\
0&0\end{smallmatrix}\right).
\]

\end{romanenumerate}

Therefore, $\Phi$ is an algebra isomorphism. Besides, it is clear
that $\Phi(\bar
X)\,\Bigl(=\Phi\bigl(-\tau(X)\bigr)\Bigr)\,=\overline{\Phi(X)}$ for
any $X\in\calT_{(\bar 1,\bar 0)}$, so that $\Phi$ is an isomorphism
of algebras with involution. This finishes the proof of the Theorem.

\bigskip

\begin{remark}\label{re:TQJ}
Given the split Cayley algebra $C$ and the unital Jordan algebra $J$
as above, the subspace $Q=\espan{1=e_1+e_2,u_0+v_0,u_1+v_1,u_2+v_2}$
is a quaternion subalgebra of $C$, which is invariant under the
action of the symmetric group $S_4$. Then so is
$\calT(Q,J)=D_{Q,Q}\oplus \bigl(Q^0\otimes J^0\bigr)\oplus d_{J,J}$,
which is a subalgebra of $\calT(C,J)$. Besides,
\[
\calT(Q,J)_{(\bar 1,\bar
0)}=kD_{u_1+v_1,u_2+v_2}\oplus\bigl((u_0+v_0)\otimes J^0\bigr).
\]
Note that $D_{u_1,u_2}=-\frac{1}{2}D_{e_2-e_1,v_0}$ and
$D_{v_1,v_2}=-\frac{1}{2}D_{e_1-e_2,u_0}$ by \eqref{eq:Du1u2}, so
that
\[
D_{u_1+v_1,u_2+v_2}=D_{v_1,u_2}+D_{u_1,v_2}-
 \frac{1}{2}D_{e_1-e_2,u_0}-\frac{1}{2}D_{e_2-e_1,v_0}.
\]
Therefore, under the isomorphism $\Phi$, $\calT(Q,J)_{(\bar 1,\bar
0)}$ maps onto the following commutative subalgebra of $\calA(J)$:
\[
S=\left\{\begin{pmatrix} 3\alpha &-\alpha +x\\
  -\alpha +x& 3\alpha\end{pmatrix}: \alpha\in k, x\in J^0\right\}.
\]
Moreover, the linear isomorphism
\[
\begin{split}
J&\longrightarrow S\\
\alpha 1+x&\mapsto \begin{pmatrix} \frac{3}{4}\alpha&
-\frac{1}{4}\alpha +\frac{1}{2}x\\
 -\frac{1}{4}\alpha +\frac{1}{2}x& \frac{3}{4}\alpha\end{pmatrix}
\end{split}
\]
for $\alpha\in k$ and $x\in J^0$, is easily checked, using
\eqref{eq:severalcross}, to be an algebra isomorphism. Hence, the
structurable algebra attached to $\calT(Q,J)$ is just, up to
isomorphism, the Jordan algebra $J$ itself.

Consider the linear isomorphism
\[
\begin{split}
\calT(Q,J)&\rightarrow \bigl(Q^0\otimes J\bigr)\oplus d_{J,J}\\
D_{a,b}&\mapsto [a,b]\otimes 1\\
a\otimes x&\mapsto a\otimes x\\
d&\mapsto d
\end{split}
\]
for any $a,b\in Q^0$, $x\in J^0$ and $d\in d_{J,J}$. Note that for
any $a,b\in Q^0$ and $x,y\in J^0$, one has
\[
[a\otimes x,b\otimes y]=t_J(x,y)D_{a,b}+[a,b]\otimes
x*y+2t(ab)d_{x,y},
\]
which maps into
\[
[a,b]\otimes(t_J(xy)+x*y)+2t(ab)d_{x,y}=([a,b]\otimes
xy)+2t(ab)d_{x,y}.
\]
Hence this linear isomorphism is an isomorphism of Lie algebras,
where the Lie bracket on $(Q^0\otimes J)\oplus d_{J,J}$ is
determined by $[a\otimes x,b\otimes y]=([a,b]\otimes
xy)+2t(ab)d_{x,y}$ for any $a,b\in Q^0$ and $x,y\in J$. This bracket
makes sense for any quaternion algebra  and any Jordan algebra (not
necessarily unital nor endowed with a normalized trace), as shown in
\cite{Tits62}. (See also Remark \ref{re:TitsTKK}.)
\end{remark}

\section{Superalgebras}

As considered in \cite{BZ96,BE03}, the Jordan algebra $J$ in Tits
construction can be replaced by a Jordan superalgebra, as long as
the superalgebra version of \eqref{eq:TCJLie} is fulfilled (see
\cite[Proposition 1.5]{BE03}). In this case $\calT(C,J)$ becomes a
Lie superalgebra. In particular \cite[Theorem 2.5]{BE03}, this is
always the case for the Jordan superalgebra $D_2$ and the Jordan
superalgebra $J=J(V,\vartheta)$ of a nondegenerate supersymmetric
superform $\vartheta$ on the superspace $V=V\subo\oplus V\subuno$
with $V\subo=0$ and $\dim V\subuno=2$. Both Jordan superalgebras are
endowed with a normalized trace.

With $C$ the split Cayley algebra over $k$,
$\calT(C,J(V,\vartheta))$ is the simple Lie superalgebra $G(3)$,
while $\calT(C,D_2)$ is the simple Lie superalgebra $F(4)$. Hence
the symmetric group $S_4$ acts on the Lie superalgebras $G(3)$ and
$F(4)$ by automorphisms. In both cases, the superalgebra with
superinvolution $\bigl(\calA(J),-\bigr)$ can be defined as in
Section \ref{se:AJ}.

The arguments in the previous Section are valid in the superalgebra
setting, as long as the convenient parity signs are inserted.
Therefore, as a consequence of Theorem \ref{th:T10AJ} and Corollary
\ref{co:AJstructurable}, we get:

\begin{theorem}\label{th:T10super}
Let $C$ be the split Cayley algebra over $k$ and let $J$ be one of
the Jordan superalgebras $J=D_2$ or $J=J(V,\vartheta)$. Let
$\calT=\calT(C,J)$ be the Lie superalgebra constructed by the Tits
construction. Then the algebra with involution $\bigl(\calT_{(\bar
1,\bar 0)},\cdot,-\bigr)$ is isomorphic to the algebra
$\bigl(\calA(J),-\bigr)$.
\end{theorem}

\begin{corollary}\label{co:TCAJsuper}
The superalgebras with superinvolution $\bigl(\calA(J),-\bigr)$, for
$J=D_2$ and $J=J(V,\vartheta)$, are structurable superalgebras.
\end{corollary}

\smallskip

\begin{remark} The simple Lie superalgebras $D(2,1;\alpha)$
($\alpha\ne 0,-1$, can be constructed directly from the Jordan
superalgebras of type $D_\alpha$ as $\calT(Q,D_\alpha)$ as in Remark
\ref{re:TQJ}, even though $D_\alpha$ has a normalized trace only for
$\alpha=2$ or $\frac{1}{2}$ (see \cite{BE03}).\quad\qed
\end{remark}

\smallskip

\begin{remark}
Consider the `tiny' Kaplansky superalgebra $K$, with even part
$K\subo=ke$, odd part $K\subuno=kx+ky$, and supercommutative
multiplication given by $e^2=e$, $ex=\frac{1}{2}x$,
$ey=\frac{1}{2}y$, and $xy=e$. $K$ is a simple nonunital Jordan
superalgebra. Then $K$ is an admissible cubic superalgebra with
$N(z)=\langle z\vert z^2\rangle$, where $\langle .\vert .\rangle$ is
the supersymmetric bilinear form such that $\langle e\vert e\rangle
=1$ and $\langle x\vert y\rangle =2$. Thus we obtain a structurable
superalgebra $\calA(K)$ as in Remark \ref{re:cubicadmissible}, with
the product given by \eqref{eq:TNNx}.

On the other hand, for the Jordan superalgebra $J=J(V,\vartheta)$,
one has $J\subo=k1$, $J\subuno=ku+kv$, with $uv=1$. Then a
straightforward computation shows that the linear map:
\[
\begin{split}
\calA(K)&\longrightarrow \calA(J)\\
\begin{pmatrix} \alpha_1&\hspace{-12pt}\gamma_1 e+\mu_1x+\nu_1y\\
  \gamma_2 e+\mu_2x+\nu_2y&\alpha_2\ \end{pmatrix}&\mapsto
  \begin{pmatrix} \alpha_1&\hspace{-12pt}\gamma_1 1-\mu_1u+2\nu_1v\\
  \gamma_2 1+\mu_2u-2\nu_2v&\alpha_2\ \end{pmatrix}
\end{split}
\]
where $\alpha_i,\gamma_i,\mu_i,\nu_i\in k$ ($i=1,2$), is an
isomorphism of structurable algebras.
\end{remark}

\section{Another action of $S_4$ on the Tits construction}

In this section, an action of the symmetric group $S_4$ by
automorphisms on the Jordan algebra of hermitian $3\times 3$
matrices over a unital composition algebra will be considered. This
is extended naturally to an action of $S_4$ by automorphisms on the
Tits construction, which gives rise to a structurable algebra. This
latter algebra is isomorphic to the tensor product of the two
composition algebras involved in the Tits construction. Therefore,
this $S_4$-action clarifies the relationship between the Tits
construction and the construction of the Lie algebras in the Magic
Square by means of a couple of composition algebras (see
\cite{BartonSudbery}, \cite{LandsbergManivel}, or \cite{Eld04}).

\smallskip

Let $\hat C$ be a unital composition algebra over our ground field
$k$, and consider the Jordan algebra of $3\times 3$ hermitian
matrices over $\hat C$:
\[
J=H_3(\hat C)=\begin{pmatrix} \alpha_1&x_0&\bar x_2\\ \bar
x_0&\alpha_2&x_1\\ x_2&\bar x_1&\alpha_0\end{pmatrix},
\]
with $\alpha_i\in k$ and $x_i\in \hat C$, $i=0,1,2$. The following
notations will be used:
\begin{gather*}
e_0=\begin{pmatrix} 0&0&0\\ 0&0&0\\ 0&0&1\end{pmatrix},\quad
e_1=\begin{pmatrix} 1&0&0\\ 0&0&0\\ 0&0&0\end{pmatrix},\quad
e_2=\begin{pmatrix} 0&0&0\\ 0&1&0\\ 0&0&0\end{pmatrix},\\
\iota_0(x)=\begin{pmatrix}0&x&0\\\bar x&0&0\\
                    0&0&0\end{pmatrix},\quad
\iota_1(x)=\begin{pmatrix}0&0&0\\ 0&0&x\\
           0&\bar x&0\end{pmatrix},\quad
\iota_2(x)=\begin{pmatrix}0&0&\bar x\\ 0&0&0\\
             x&0&0\end{pmatrix}.
\end{gather*}
$J$ is a Jordan algebra with the product given by $X\circ
Y=\frac{1}{2}(XY+YX)$, which satisfies:
\[
\begin{split}
&e_i\circ e_j=\delta_{ij}e_i,\\
&e_{i+1}\circ \iota_i(x)=e_{i+2}\circ
\iota_i(x)=\frac{1}{2}\iota_i(x),\quad
   e_i\circ\iota_i(x)=0,\\
&\iota_i(x)\circ\iota_{i+1}(y)=\frac{1}{2}\iota_{i+2}(\overline{xy}),\\
&\iota_i(x)\circ \iota_i(y)=\frac{1}{2}t(x\bar y)(e_{i+1}+e_{i+2}),
\end{split}
\]
where $x,y\in \hat C$, $t$ denotes the trace in $\hat C$,
$x\mapsto\bar x$ the canonical involution, and the indices are taken
modulo $3$.

The symmetric group $S_4$ embeds in the automorphism group of $J$ as
follows:
\[
\begin{split}
\tau_1&: e_i\mapsto e_i,\ \iota_0(x)\mapsto \iota_0(x),\
  \iota_1(x)\mapsto -\iota_1(x),\ \iota_2(x)\mapsto
  -\iota_2(x),\\[2pt]
\tau_2&: e_i\mapsto e_i,\ \iota_0(x)\mapsto -\iota_0(x),\
 \iota_1(x)\mapsto \iota_1(x),\ \iota_2(x)\mapsto
 -\iota_2(x),\\[2pt]
\varphi&: e_i\mapsto e_{i+1},\ \iota_i(x)\mapsto
\iota_{i+1}(x),\\[2pt]
\tau&: e_0\mapsto e_0,\ e_1\mapsto e_2,\ e_2\mapsto e_1,\
\iota_0(x)\mapsto \iota_0(\bar x),\\
 &\hspace{140pt} \iota_1(x)\mapsto \iota_2(\bar x),\ \iota_2(x)\mapsto \iota_1(\bar x),
\end{split}
\]
for any $i=0,1,2$ (indices modulo $3$) and $x\in\hat C$.

\smallskip

The arguments in \cite[Chapter IV, \S 9]{Schafer} show that there is
the following grading over $\bZ_2\times\bZ_2$ of the Lie algebra of
derivations of $J$:
\[
\der J=d_{J,J}=\{d\in\der J: d(e_i)=0\ (i=0,1,2)\}\oplus
\bigl(\oplus_{i=0}^2d_{e_{i+1}-e_{i+2},\iota_i(\hat C)}\bigr),
\]
where $d_{x,y}$ is defined as in \eqref{eq:dxy}: $d_{x,y}(z)=x\circ
(y\circ z)-y\circ(x\circ z)$. This is precisely the grading induced
by the action of Klein's $4$-group.

\smallskip

The action of $S_4$ on $J$ extends to an action of $S_4$ by
automorphisms of the Lie algebra $\calT(C,J)$, where $C$ is another
unital composition algebra over $k$:
\begin{equation}\label{eq:S4actiononTCJbis}
\psi\bigl(D+(a\otimes x)+d\bigr)= D+(a\otimes \psi(x))+\psi
d\psi^{-1},
\end{equation}
for any $\psi\in S_4$, $D\in \der C$, $d\in d_{J,J}$, $a\in C^0$ and
$x\in J^0$.

As in \eqref{eq:CZ2Z2grading} and \eqref{eq:TZ2Z2grading}, the
action of Klein's $4$-group induces a grading over $\bZ_2\otimes
\bZ_2$:
\begin{equation}\label{eq:TZ2Z2gradingbis}
\calT=\calT(C,J)=\calT_{(\bar 0,\bar 0)}\oplus \calT_{(\bar 1,\bar
0)}\oplus \calT_{(\bar 0,\bar 1)}\oplus \calT_{(\bar 1,\bar 1)}.
\end{equation}
Again (see \cite[Section 2]{EOJAlgebra}), the subspace $\calT_{(\bar
1,\bar 0)}=\{ X\in\calT: \tau_1(X)=X,\ \tau_2(X)=-X\}$ becomes an
algebra with involution by means of:
\begin{equation}\label{eq:T10productbis}
\begin{split}
&X\cdot Y=-\tau\bigl([\varphi(X),\varphi^2(Y)]\bigr),\\
&\bar X=-\tau(X),
\end{split}
\end{equation}
for any $X,Y\in\calT_{(\bar 1,\bar 0)}$, as in
\eqref{eq:T10product}.

\smallskip

\begin{theorem}\label{th:T10CChat}
Let $C$ and $\hat C$ be two unital composition algebras over $k$,
whose traces will be both denoted by $t$, and let $J$ be the Jordan
algebra of $3\times 3$ hermitian matrices over $\hat C$. Let
$\calT=\calT(C,J)$ be the associated Tits algebra, with the action
of the symmetric group $S_4$ given by \eqref{eq:S4actiononTCJbis}.
Then the algebra with involution $\bigl(\calT_{(\bar 1,\bar
0)},\cdot,-\bigr)$ is isomorphic to the structurable algebra
$C\otimes \hat C$, with multiplication $(a\otimes x)(b\otimes
y)=ab\otimes xy$ and involution $\overline{a\otimes x}=\bar
a\otimes\bar x$, for any $a,b\in C$ and $x,y\in\hat C$.
\end{theorem}
\begin{proof}
To begin with, the subspace $\calT_{(\bar 1,\bar 0)}$ is
\[
\calT_{(\bar 1,\bar 0)}=\bigl(C^0\otimes \iota_0(\hat C)\bigr)\oplus
d_{e_1-e_2,\iota_0(\hat C)}.
\]
For ease of notation, write $d_i(x)=d_{e_{i+1}-e_{i+2},\iota_i(x)}$
for any $i=0,1,2$ (modulo $3$) and $x\in \hat C$, so
\[
\calT_{(\bar 1,\bar 0)}=\bigl(C^0\otimes \iota_0(\hat C)\bigr)\oplus
d_0(\hat C).
\]
Consider the bijective linear map
\[
\begin{split}
\Phi:\calT_{(\bar 1,\bar 0)}&\rightarrow C\otimes\hat C\\
 a\otimes \iota_0(x)&\mapsto -a\otimes x,\\
 d_0(x)&\mapsto -\frac{1}{2}1\otimes x,
\end{split}
\]
for $a\in C^0$ and $x\in\hat C$.

It is clear the $\Phi(\bar X)=\overline{\Phi(X)}$ for any
$X\in\calT_{(\bar 1,\bar 0)}$ since
\[
\begin{split}
&-\tau(a\otimes \iota_0(x))=-a\otimes \iota_0(\bar x)=\bar a\otimes
\iota_0(\bar x),\ \text{and}\\
 &-\tau(d_0(x))=-[l_{\tau(e_1-e_2)},l_{\tau(\iota_0(x))}]=
  -[l_{e_2-e_1},l_{\iota_0(\bar x)}]=d_0(\bar x),
\end{split}
\]
for any $a\in C^0$ and $x\in \hat C$. Note that the standard
involutions of both $C$ and $\hat C$ are denoted by the same symbol.

To prove that $\Phi$ is an algebra homomorphism, the following
instances of $\Phi(X\cdot Y)=\Phi(X)\Phi(Y)$ have to be checked:

\smallskip

\begin{romanenumerate}

\item $X=d_0(x)$, $Y=d_0(y)$, with $x,y\in C$:\quad Note that one
has
\[
\varphi(d_0(x))=\varphi\bigl([l_{e_1-e_2},l_{\iota_0(x)}]\bigr)=
 [l_{\varphi(e_1-e_2)},l_{\varphi(\iota_0(x))}]=
 [l_{e_2-e_0},l_{\iota_1(x)}]=d_1(x),
\]
while $\varphi^2(d_0(y))=d_2(y)$. Also, the equality
$\tau(d_0(x))=-d_0(\bar x)$ was checked above. Hence,
\[
X\cdot Y=-\tau\bigl([\varphi(d_0(x)),\varphi^2(d_0(y))]\bigr)=
 -\tau\bigl([d_1(x),d_2(y)]\bigr).
\]
But $[d_1(x),d_2(y)]$ belongs to the subspace $d_0(\hat C)$ (because
of the grading of $\der J$ over $\bZ_2\times \bZ_2$, so there is an
element $z\in\hat C$ such that  $[d_1(x),d_2(y)]=d_0(z)$. Now, for
any $z\in\hat C$ and any $i=0,1,2$, one has:
\[
\begin{split}
d_i(z)(e_{i+1})&=[l_{e_{i+1}-e_{i+2}},l_{\iota_i(z)}](e_{i+1})\\
 &=(e_{i+1}-e_{i+2})\circ(\iota_i(z)\circ e_{i+1})-
   \iota_i(z)\circ((e_{i+1}-e_{i+2})\circ e_{i+1})\\
 &=\frac{1}{2}(e_{i+1}-e_{i+2})\circ\iota_i(z)
   -\iota_i(z)\circ e_{i+1}\\
 &=-\frac{1}{2}\iota_i(z).
\end{split}
\]
In the same vein, one obtains:
\begin{equation}\label{eq:dizej}
d_i(z)(e_{i+1})=-d_i(z)(e_{i+2})=-\frac{1}{2}\iota_i(z),\quad
d_i(z)(e_i)=0,
\end{equation}
for any $z\in\hat C$ and $i=0,1,2$. Also,
\[
\begin{split}
d_i(z)&(\iota_{i+1}(t))\\
&=
 (e_{i+1}-e_{i+2})\circ(\iota_i(z)\circ\iota_{i+1}(t))-
    \iota_i(z)\circ((e_{i+1}-e_{i+2})\circ\iota_{i+1}(t))\\
 &=\frac{1}{2}(e_{i+1}-e_{i+2})\circ\iota_{i+2}(\overline{zt})+
   \frac{1}{2}\iota_i(z)\circ\iota_{i+1}(t)\\
 &=\frac{1}{4}\iota_{i+2}(\overline{zt})
    +\frac{1}{4}\iota_{i+2}(\overline{zt})=\frac{1}{2}\iota_{i+2}(\overline{zt}),
\end{split}
\]
and, in the same vein,
\begin{equation}\label{eq:dizdi1t}
\begin{split}
d_i(z)(\iota_{i+1}(t))&=\frac{1}{2}\iota_{i+2}(\overline{zt}),\\
d_i(z)(\iota_{i+2}(t))&=-\frac{1}{2}\iota_{i+1}(\overline{tz}),\\
d_i(z)(\iota_{i}(t))&=\frac{1}{2}t(z\bar t)(e_{i+1}-e_{i+2})
\end{split}
\end{equation}
for any $z,t\in\hat C$ and $i=0,1,2$.

Thus,
\[
\begin{split}
[d_1(x),d_2(y)](e_1)&=d_1(x)\bigl(d_2(y)(e_1)\bigr)-
   d_2(y)\bigl(d_1(x)(e_1)\bigr)\\
   &=d_1(x)\bigl(\frac{1}{2}\iota_2(y)\bigr)=
   \frac{1}{4}\iota_0(\overline{xy}),
\end{split}
\]
and therefore,
\[
[d_1(x),d_2(y)]=-\frac{1}{2}d_0(\overline{xy})
\]
because of \eqref{eq:dizej}. Hence
\[
X\cdot
Y=-\tau\big([d_1(x),d_2(y)]\bigr)
 =\tau\bigl(\frac{1}{2}d_0(\overline{xy})\bigr)
 =-\frac{1}{2}d_0(xy),
\]
and hence
\[
\Phi(X\cdot Y)=\frac{1}{4}1\otimes xy=\bigl(-\frac{1}{2}1\otimes
x\bigr)\bigl(-\frac{1}{2} 1\otimes y\bigr)=\Phi(X)\Phi(Y).
\]

\medskip

\item $X=d_0(x)$, $Y=a\otimes \iota_0(y)$, $x,y\in\hat C$, $a\in
C^0$.\quad Here
\[
\begin{split}
X\cdot Y&=
-\tau\bigl([\varphi(X),\varphi^2(Y)]\bigr)
   =-\tau\bigl([d_1(x),a\otimes\iota_2(y)]\bigr)\\
   &=-\tau\bigl(a\otimes d_1(x)(\iota_2(y))\bigr)=
     -\tau\bigl(a\otimes\frac{1}{2}\iota_0(\overline{xy})\bigr)\quad\text{(by
     \eqref{eq:dizdi1t})}\\
   &=-\frac{1}{2}a\otimes\iota_0(xy),
\end{split}
\]
so
\[
\Phi(X\cdot Y)=\frac{1}{2}a\otimes xy=\bigl(-\frac{1}{2}1\otimes
x\bigr)\bigl(-a\otimes y)=\Phi(X)\Phi(Y).
\]

\medskip

\item $X=a\otimes \iota_0(x)$, $Y=d_0(y)$, $a\in C^0$, $x,y\in \hat
C$.\quad In this case,
\[
\begin{split}
X\cdot Y&=-\tau\bigl([\varphi(X),\varphi^2(Y)]\bigr)=
  -\tau\bigl([a\otimes\iota_1(x),d_2(y)]\bigr)\\
  &=\tau\bigl(a\otimes d_2(y)(\iota_1(x))\bigr)=
   \tau\bigl(-\frac{1}{2}a\otimes\iota_0(\overline{xy})\bigr)\quad\text{(by
   \eqref{eq:dizdi1t})}\\
  &=-\frac{1}{2} a\otimes\iota_0(xy),
\end{split}
\]
so
\[
\Phi(X\cdot Y)=\frac{1}{2}a\otimes xy=\bigl(-a\otimes
x\bigr)\bigl(-\frac{1}{2}1\otimes y\bigr)=\Phi(X)\Phi(Y).
\]

\medskip

\item $X=a\otimes \iota_0(x)$, $Y=b\otimes \iota_0(y)$, for $a,b\in
C^0$ and $x,y\in \hat C$.\quad Here
\[
\begin{split}
X\cdot Y&= -\tau\bigl([\varphi(X),\varphi^2(Y)]\bigr)=
  -\tau\bigl([a\otimes\iota_1(x),b\otimes\iota_2(y)]\bigr)\\
  &=-\tau\bigl(([a,b]\otimes\frac{1}{2}\iota_0(\overline{xy}))\, +\,
   2t(ab)[l_{\iota_1(x)},l_{\iota_2(y)}]\bigr).
\end{split}
\]
But $[l_{\iota_1(x)},l_{\iota_2(y)}]=d_0(z)$ for some $z\in\hat C$,
and $d_0(z)(e_1)=-\frac{1}{2}\iota_0(z)$ by \eqref{eq:dizej}, while
\[
\begin{split}
[l_{\iota_1(x)},l_{\iota_2(y)}](e_1)&=
   \iota_1(x)\circ(\iota_2(y)\circ
   e_1)-\iota_2(y)\circ(\iota_1(x)\circ e_1)\\
   &=\frac{1}{2}\iota_1(x)\circ\iota_2(y)=\frac{1}{4}\iota_0(\overline{xy}).
\end{split}
\]
Hence
$[l_{\iota_1(x)},l_{\iota_2(y)}]=-\frac{1}{2}d_0(\overline{xy})$,
and
\[
\begin{split}
X\cdot Y&=
-\tau\bigl(([a,b]\otimes\frac{1}{2}\iota_0(\overline{xy}))\, -\,
t(ab)d_0(\overline{xy})\bigr)\\
 &=-\frac{1}{2}([a,b]\otimes \iota_0(xy))\,-\, t(ab)d_0(xy).
\end{split}
\]
Therefore,
\[
\begin{split}
\Phi(X\cdot Y)&=\frac{1}{2}([a,b]\otimes
xy)+\frac{1}{2}t(ab)(1\otimes xy)\\
&=\frac{1}{2}([a,b]+t(ab)1)\otimes xy.
\end{split}
\]
But $ab+ba=t(ab)1$ for any $a,b\in C^0$, while $ab-ba=[a,b]$. Hence
$ab=\frac{1}{2}([a,b]+t(ab)1)$, and
\[
\Phi(X\cdot Y)=ab\otimes xy=(-a\otimes x)(-b\otimes
y)=\Phi(X)\Phi(Y)
\]
also in this case.
\end{romanenumerate}
\end{proof}

\section{$S_4$-actions and structurable algebras}

Among the irreducible representations of the symmetric group $S_4$,
let us consider the one obtained on the tensor product of the
standard representation and the alternating one \cite[\S
2.3]{FultonHarris}. This is obtained on a three dimensional vector
space $W=kw_0+kw_1+kw_2$ with the action of $S_4$ given by
\begin{equation}\label{eq:S4W}
\left\{\begin{aligned} \tau_1&: w_0\mapsto w_0,\ w_1\mapsto -w_1,\
w_2\mapsto -w_2,\\
\tau_2&: w_0\mapsto -w_0,\ w_1\mapsto w_1,\ w_2\mapsto -w_2,\\
\varphi&: w_0\mapsto w_1\mapsto w_2\mapsto w_0,\\
\tau&: w_0\mapsto -w_0,\ w_1\mapsto -w_2,\ w_2\mapsto -w_1.
\end{aligned}\right.
\end{equation}
(This is the representation that appears on the subspaces spanned by
the $u_i$'s and the $v_i$'s in \eqref{eq:S4actiononC}.)

Then the general Lie algebra $\frgl(W)$ becomes a module for $S_4$:
$\sigma\cdot f=\sigma f\sigma^{-1}$. Thus, $S_4$ acts by
automorphisms on $\frgl(W)$. Identifying $\frgl(W)$ with $\Mat_3(k)$
by means of our basis $\{w_1,w_2,w_0\}$, consider the following
basis of $\frgl(W)$:
\begin{equation}\label{eq:HGDs}
\begin{aligned}
H_0&=\begin{pmatrix} 0&0&0\\ 0&0&0\\ 0&0&1\end{pmatrix}&
 H_1&=\begin{pmatrix} 1&0&0\\ 0&0&0\\ 0&0&0\end{pmatrix}&
 H_2&=\begin{pmatrix} 0&0&0\\ 0&1&0\\ 0&0&0\end{pmatrix}\\
G_0&=\begin{pmatrix} 0&1&0\\ 1&0&0\\ 0&0&0\end{pmatrix}&
 G_1&=\begin{pmatrix} 0&0&0\\ 0&0&1\\ 0&1&0\end{pmatrix}&
 G_2&=\begin{pmatrix} 0&0&1\\ 0&0&0\\ 1&0&0\end{pmatrix}\\
 D_0&=\begin{pmatrix} 0&-1&0\\ 1&0&0\\ 0&0&0\end{pmatrix}&
 D_1&=\begin{pmatrix} 0&0&0\\ 0&0&-1\\ 1&0&0\end{pmatrix}&
 D_2&=\begin{pmatrix} 0&0&1\\ 0&0&0\\ -1&0&0\end{pmatrix}
 \end{aligned}
\end{equation}

The space $W$ is endowed with a natural nondegenerate symmetric
bilinear form $(.\vert .):W\times W\rightarrow k$ given by
$(w_i\vert w_j)=\delta_{ij}$, which is invariant under the action of
$S_4$. Actually, $S_4$ embeds in the associated special orthogonal
group $SO(W)$. The corresponding orthogonal Lie algebra $\frso_3=\{
F\in\frgl(W): (F(v)\vert w)+(v\vert F(w))=0\ \forall v,w\in W\}$ is
the span of the $D_i$'s in \eqref{eq:HGDs}. Note that
\[
[D_i,D_{i+1}]=D_{i+2}\quad \text{(indices modulo $3$).}
\]

As a module for $\frso_3$, $\frgl(W)$ decomposes into the following
direct sum of irreducible modules (remember that the characteristic
of the ground field $k$ is assumed to be $\ne 2,3$):
\begin{equation}\label{eq:glWso3hz}
\frgl(W)=\frso_3 \oplus \frh \oplus \frz,
\end{equation}
where $\frz=kI_3$ ($I_3$ denotes the identity matrix), and $\frh=\{
F\in\frgl(W):(F(v)\vert w)=(v\vert F(w))\ \forall v,w\in W\
\text{and} \trace(F)=0\}$.

These three irreducible modules: $\frso_3$, $\frh$, and $\frz$, are
invariant under the action by conjugation by the orthogonal group,
and hence, in particular, under the action of $S_4$, but while
$\frso_3$ and $\frz$ are irreducible modules under the action of
$S_4$, $\frh$ decomposes as the direct sum of two irreducible
modules for $S_4$:
\[
\frh=\espan{G_0,G_1,G_2}\oplus\espan{H_0-H_1,H_1-H_2}.
\]
A simple computation shows that $\espan{H_0,H_1,H_2}$ is left
elementwise fixed by Klein's $4$-group $V$, and becomes the natural
module for $S_3=S_4/V$. On the other hand, $\espan{G_0,G_1,G_2}$ is
the standard module for $S_4$. Thus, among the five irreducible
modules for $S_4$ (up to isomorphism), only the alternating one is
missing in $\frgl(W)$.

\begin{lemma}\label{le:so3invariantmaps}
Up to scalars, the following maps are the unique $\frso_3$-invariant
linear maps between the $\frso_3$-modules considered:

\begin{itemize}
\item $\frso_3\otimes\frso_3\rightarrow \frso_3:$ $A\otimes
B\mapsto [A,B]$,

\smallskip

\item $\frso_3\otimes\frso_3\rightarrow \frh:$ $A\otimes B\mapsto
AB+BA-\frac{2}{3}\trace(AB)I_3$,

\smallskip

\item $\frso_3\otimes\frso_3\rightarrow \frz:$ $A\otimes B\mapsto
\trace(AB)I_3$,

\smallskip

\item $\frso_3\otimes \frh\rightarrow \frso_3:$ $A\otimes X\mapsto
AX+XA$,

\smallskip

\item $\frso_3\otimes \frh\rightarrow \frh:$ $A\otimes X\mapsto
[A,X]$,

\smallskip

\item $\frso_3\otimes \frh\rightarrow \frz:$ $A\otimes X\mapsto 0$,

\smallskip

\item $\frh\otimes \frh\rightarrow \frso_3:$ $X\otimes Y\mapsto
[X,Y]$,

\smallskip

\item $\frh\otimes \frh\rightarrow \frh:$ $X\otimes Y\mapsto
XY+YX-\frac{2}{3}\trace(XY)I_3$,

\smallskip

\item $\frh\otimes \frh\rightarrow\frz:$ $X\otimes Y\mapsto
\trace(XY)I_3$.
\end{itemize}

Moreover, all these maps are invariant under the action of $S_4$.
\end{lemma}

\begin{proof}
It is clear that all these maps are invariant under the action of
both $\frso_3$ and $S_4$, because so is the trace form and the
associative multiplication in $\End_k(W)$.

Now, to prove the uniqueness it is enough to assume the ground field
$k$ to be algebraically closed. In this case, $\frso_3$ is
isomorphic to $\frsl_2$ and, as a module for $\frsl_2$, $\frso_3$ is
isomorphic to $V(2)$, $\frh$ to $V(4)$ and $\frz$ to $V(0)$
(notation as in \cite[\S 7]{Humphreys}). Note that this makes sense
because the characteristic is either $0$ or $\geq 5$. But for
$n,m=0,2$ or $4$, $V(m)\otimes V(n)$ is generated, as a module for
$\frsl_2$, by $R\otimes S$, for a highest weight vector $R$ of
$V(m)$ and a lowest weight vector $S$ of $V(n)$. Hence any invariant
linear map $V(m)\otimes V(n)\rightarrow V(p)$ ($p=0,2$, or $4$) is
determined by the image of $R\otimes S$, which belongs to the weight
space of $V(p)$ of weight $2(m-n)$. This is at most one-dimensional,
and the result follows.
\end{proof}

\medskip

Let $\frg$ be an arbitrary Lie algebra over $k$ endowed with an
action of the symmetric group $S_4$ by automorphisms, that is,
endowed with a group homomorphism
\[
S_4\rightarrow \Aut(\frg).
\]
As before, the action of Klein's $4$-group induces a grading of
$\frg$ over $\bZ_2\times \bZ_2$ and the subspace $A=\frg_{(\bar
1,\bar 0)}$ becomes an algebra with involution by means of:
\[
\left\{\begin{aligned} &x\cdot
y=-\tau\bigl([\varphi(x),\varphi^2(y)]\bigr),\\[2pt]
&\bar x=-\tau(x).\end{aligned}\right.
\]
The algebra $(A,\cdot,-)$ will be called the \emph{coordinate
algebra} of $\frg$. If this algebra is unital, then it is
structurable (\cite[Theorem 2.6]{OkuboLAA} and \cite[Theorem
2.9]{EOJAlgebra}). This is the situation that has already appeared
in Theorems \ref{th:T10AJ} and \ref{th:T10CChat}.

\begin{theorem}\label{th:g135}
Let $\frg$ be a Lie algebra over $k$. Then $\frg$ is endowed with an
action of $S_4$ by automorphisms such that the coordinate algebra is
unital (or, equivalently, structurable) if and only if there is a
subalgebra of $\frg$ isomorphic to $\frso_3$, such that, as a module
for this subalgebra, $\frg$ is the direct sum of irreducible modules
isomorphic either to the adjoint module $\frso_3$, the five
dimensional module $\frh$ or the trivial one dimensional module
$\frz$.
\end{theorem}

\begin{proof}
Assume first that $\frg$ contains $\frso_3$ as a subalgebra with the
properties stated in the Theorem. Then, collecting isomorphic
irreducible modules, we may write:
\begin{equation}\label{eq:gso3hd}
\frg= \bigl(\frso_3\otimes\calH\bigr)\oplus\bigl(\frh\otimes
\calS\bigr)\oplus\frd,
\end{equation}
for vector subspaces $\calH$, $\calS$ and $\frd$. The subalgebra
$\frso_3$ is then identified to $\frso_3\otimes 1$ for a
distinguished element $1\in\calH$. Here $\frd=\{x\in \frg: [d,x]=0\
\forall d\in\frso_3\}$ is the sum of the trivial irreducible
modules, so $\frd$ is the centralizer of the subalgebra $\frso_3$
and, in particular, it is a subalgebra of $\frg$.

Because of Lemma \ref{le:so3invariantmaps}, the Lie bracket in
$\frg$, which is invariant under the action of the subalgebra
$\frso_3$, is given by:

\begin{itemize}

\item $\frd$ is a subalgebra of $\frg$,

\smallskip

\item $[A\otimes a,B\otimes b]=[A,B]\otimes a\circ b\, -\,
\bigl(AB+BA-\frac{2}{3}\trace(AB)I_3\bigr)\otimes\frac{1}{2}[a,b]\,
+\trace(AB)d_{a,b}$,

\smallskip

\item $[A\otimes a,X\otimes x]=-(AX+XA)\otimes \frac{1}{2}[a,x]\,
+\, [A,X]\otimes a\circ x$,

\smallskip

\item $[X\otimes x,Y\otimes y]=[X,Y]\otimes x\circ y\
-\bigl(XY+YX-\frac{2}{3}\trace(XY)I_3\bigr)\otimes
\frac{1}{2}[x,y]\, +\trace(XY)d_{x,y}$,

\smallskip

\item $[d,A\otimes a]=A\otimes d(a)$,

\smallskip

\item  $[d,X\otimes x]=X\otimes d(x)$,
\end{itemize}

\noindent for any $A,B\in\frso_3$, $X,Y\in\frh$, $a,b\in\calH$,
$x,y\in\calS$, and $d\in \frd$, where

\begin{itemize}
\item[--] $\calH\times\calH\rightarrow \calH$: $(a,b)\mapsto a\circ
b$ is a symmetric bilinear map with $1\circ a=a$ for any
$a\in\calH$,

\item[--] $\calH\times\calH\rightarrow \calS$: $(a,b)\mapsto [a,b]$
is a skew symmetric bilinear map with $[1,a]=0$ for any $a\in\calH$,

\item[--] $\calH\times\calS\rightarrow \calH$: $(a,x)\mapsto [a,x]$
is a bilinear map with $[1,x]=0$ for any $x\in\calS$,

\item[--] $\calH\times\calS\rightarrow \calS$: $(a,x)\mapsto a\circ x$
is a bilinear map with $1\circ x=x$ for any $x\in\calS$,

\item[--] $\calS\times \calS\rightarrow \calH$: $(x,y)\mapsto x\circ
y$ is a symmetric bilinear map,

\item[--] $\calS\times\calS\rightarrow \calS$: $(x,y)\mapsto [x,y]$
is a skew symmetric bilinear map,

\item[--] $\calH\times\calH\rightarrow \frd$: $(a,b)\mapsto d_{a,b}$
is a skew symmetric bilinear map,

\item[--] $\calS\times\calS\rightarrow \frd$: $(x,y)\mapsto d_{x,y}$
is a skew symmetric bilinear map,

\item[--] the bilinear maps $\frd\times\calH\rightarrow\calH$:
$(d,a)\mapsto d(a)$ and $\frd\times\calS\rightarrow \calS$:
$(d,x)\mapsto d(x)$, give two representations of the Lie algebra
$\frd$.

\end{itemize}
\smallskip

Now, define an action of $S_4$ on $\frg$ by means of the actions by
conjugation of $S_4$ on both $\frso_3$ and $\frh$:
\begin{equation}\label{eq:S4actiong}
\psi\bigl(A\otimes a + X\otimes x + D\bigr)=(\psi\cdot A)\otimes a +
(\psi\cdot X)\otimes x + D
\end{equation}
for any $\psi\in S_4$, $A\in\frso_3$, $X\in\frh$, $a\in\calH$,
$x\in\calS$ and $D\in\frd$.

The invariance of the maps in Lemma \ref{le:so3invariantmaps} under
the action of $S_4$ immediately implies that any $\psi\in S_4$ acts
as an automorphism of $\frg$.

Besides, the subspace $\frg_{(\bar 1,\bar 0)}=\{ g\in\frg :
\tau_1(g)=g=-\tau_2(g)\}$ is precisely the subspace
\[
D_0\otimes\calH\, \oplus \, G_0\otimes\calS.
\]
The involution in the coordinate algebra is given by
\[
\overline{D_0\otimes a + G_0\otimes x}=-\tau(D_0)\otimes
a-\tau(G_0)\otimes x=D_0\otimes a-G_0\otimes x,
\]
for any $a\in\calH$ and $x\in\calS$, and the multiplication in the
coordinate algebra is given by:
\[
\begin{split}
\bigl(D_0\otimes a + \null&G_0\otimes x\bigr)\cdot
 \bigl(D_0\otimes b + G_0\otimes y\bigr)\\
 &=-\tau\bigl([\varphi(D_0\otimes a+G_0\otimes
 x),\varphi^2(D_0\otimes b +G_0\otimes y)]\bigr)\\
 &=-\tau\bigl([D_1\otimes a+G_1\otimes x,D_2\otimes b+G_2\otimes
 y]\bigr).
\end{split}
\]
But,
\[
\begin{split}
&[D_1,D_2]=D_0,\ [D_1,G_2]=-G_0,\ [D_2,G_1]=G_0,\
[G_1,G_2]=D_0,\\[4pt]
&D_1D_2+D_2D_1-\frac{2}{3}\trace(D_1D_2)I_3=G_0,\ D_1G_2+G_2D_1=-D_0,\\
&D_2G_1+G_1D_2=-D_0,\
G_1G_2+G_2G_1-\frac{2}{3}\trace(G_1G_2)I_3=G_0,\\[4pt]
&\trace(D_1D_2)=0=\trace(G_1G_2),
\end{split}
\]
so
\[
\begin{split}
\bigl(D_0\otimes a + \null&G_0\otimes x\bigr)\cdot
 \bigl(D_0\otimes b + G_0\otimes y\bigr)\\
 &=-\tau\bigl([D_1\otimes a+G_1\otimes x,D_2\otimes b+G_2\otimes
 y]\bigr)\\[4pt]
 &=-\tau\bigl( D_0\otimes a\circ b -
 G_0\otimes\frac{1}{2}[a,b]+D_0\otimes\frac{1}{2}[a,y]-G_0\otimes
 a\circ y\\
 &\qquad -D_0\otimes b\circ x + G_0\otimes\frac{1}{2}[b,x] +D_0\otimes
 x\circ y - G_0\otimes\frac{1}{2}[x,y]\bigr)\\[4pt]
 &=\bigl(D_0\otimes a\circ b+G_0\otimes\frac{1}{2}[a,b]\bigr)
  +\bigl(D_0\otimes\frac{1}{2}[a,y]+G_0\otimes\frac{1}{2}a\circ
  y\bigr)\\
  &\qquad -\bigl(D_0\otimes b\circ x-G_0\otimes\frac{1}{2}[b,x]\bigr)
  +\bigl(D_0\otimes x\circ y+G_0\otimes\frac{1}{2}[x,y]\bigr).
\end{split}
\]

Define $x\circ b=b\circ x$ and $[x,b]=-[b,x]$ for any $b\in \calH$
and $x\in\calS$. Now consider the vector space
$\calA=\calH\oplus\calS$ and define a multiplication on it by means
of
\[
u\cdot v=u\circ v+\frac{1}{2}[u,v]
\]
for any $u,v\in \calH\cup\calS$, so $u\circ v=\frac{1}{2}(u\cdot
v+v\cdot u)$ and $[u,v]=u\cdot v-v\cdot u$. Define too a linear map
$-:\calA\rightarrow \calA$ such that $\overline{a+x}=a-x$ for any
$a\in\calH$ and $x\in\calS$. Then the linear map $D_0\otimes
a+G_0\otimes x\mapsto a+x$ gives an isomorphism between the
coordinate algebra $\frg_{(\bar 1,\bar 0)}$ and the algebra with
involution $(\calA,\cdot,-)$. Besides, $1\in\calH$ is the unity
element of $\calA$.

\smallskip

Conversely, let $\frg$ be a Lie algebra with an action of $S_4$ by
automorphisms such that the coordinate algebra is unital. As in
\cite{EOJAlgebra}, let $\calA=\frg_{(\bar 1,\bar 0)}$ be the
coordinate algebra, and for any $x\in\calA$ consider the elements:
\[
\iota_0(x)=x\in\frg_{(\bar 1,\bar 0)},\quad
\iota_1(x)=\varphi(x)\in\frg_{(\bar 0,\bar 1)},\quad
\iota_2(x)=\varphi^2(x)\in\frg_{(\bar 1,\bar 1)}.
\]
(Recall that $\varphi$ is the cycle $(123)$ in $S_4$.)

Then \cite[\S 2]{EOJAlgebra}, for any $x,y\in\calA$ and $i=0,1,2$
(indices modulo $3$):
\[
[\iota_i(x),\iota_{i+1}(y)]=\iota_{i+2}(\overline{x\cdot y}).
\]

Therefore, $\frs=\espan{\iota_0(1),\iota_1(1),\iota_2(1)}$ is a
subalgebra of $\frg$ isomorphic to $\frso_3$ (by means of
$\iota_i(1)\mapsto D_i$ for any $i=0,1,2$). This is the subalgebra
we are looking for.

For any $0\ne x\in\calA$ with $\bar x=x$,
$\espan{\iota_0(x),\iota_1(x),\iota_2(x)}$ is a copy of the adjoint
module for $\frs$, because
\begin{equation}\label{eq:iota1iotaxsym}
\begin{split}
[\iota_i(1),\iota_{i+1}(x)]
   &=\iota_{i+2}(x)=[\iota_i(x),\iota_{i+1}(1)],\\[4pt]
[\iota_i(1),\iota_i(x)]&=[[\iota_{i+1}(1),\iota_{i+2}(1)],\iota_i(x)]\\
 &=[[\iota_{i+1}(1),\iota_i(x)],\iota_{i+2}(1)]+
   [\iota_{i+1}(1),[\iota_{i+2}(1),\iota_{i}(x)]]\\
 &=-[\iota_{i+2}(x),\iota_{i+2}(1)]+[\iota_{i+1}(1),\iota_{i+1}(x)]\\
 &=[\iota_{i+1}(1),\iota_{i+1}(x)]+[\iota_{i+2}(1),\iota_{i+2}(x)],
\end{split}
\end{equation}
for any $i=0,1,2$. Adding the resulting equations for $i=0,1,2$
gives $\sum_{i=0}^2[\iota_i(1),\iota_i(x)]
=2\bigl(\sum_{i=0}^2[\iota_i(1),\iota_i(x)]\bigr)$, so
$\sum_{i=0}^2[\iota_i(1),\iota_i(x)]=0$ and then
\eqref{eq:iota1iotaxsym} implies that $[\iota_i(1),\iota_i(x)]=0$
for any $i$.

Now take an element $0\ne x\in\calA$ with $\bar x=-x$. Let us first
check that the $\frs$-submodule generated by $\iota_0(x)$ is
$\frv=\espan{\iota_i(x),[\iota_i(1),\iota_i(x)], i=0,1,2}$. To do
so, by symmetry, it is enough to check that this subspace is closed
under the action of $\iota_0(1)$, but:
\begin{equation}\label{eq:iota1iotaxskew}
\begin{split}
&[\iota_0(1),\iota_1(x)]=\iota_2(\bar x)=-\iota_2(x),\\
&[\iota_0(1),\iota_2(x)]=-\iota_1(\bar x)=\iota_1(x),\\[4pt]
&[\iota_0(1),[\iota_1(1),\iota_1(x)]]
   =[[\iota_0(1),\iota_1(1)],\iota_1(x)]+
   [\iota_1(1),[\iota_0(1),\iota_1(x)]]\\
&\phantom{[\iota_0(1),[\iota_1(1),\iota_1(x)]]}
   =[\iota_2(1),\iota_1(x)]-[\iota_1(1),\iota_2(x)]=2\iota_0(x),\\
&[\iota_0(1),[\iota_2(1),\iota_2(x)]]=2\iota_0(x)\ \text{(same
arguments)},
\end{split}
\end{equation}
and finally, as in \cite[Theorem 2.4]{EOJAlgebra},
$[\iota_0(1),[\iota_0(1),\iota_0(x)]]=-\iota_0(\delta_0(1,x)(1))$,
with $\delta_0(1,x)=-\delta_1(\bar x,1)-\delta_2(1,x)=2(L_x+R_x)$,
where $L_x$ and $R_x$ denote, respectively, the left and right
multiplication by $x$ in $\calA$. Hence,
\begin{equation}\label{eq:iota1iotaxskewbis}
[\iota_0(1),[\iota_0(1),\iota_0(x)]]=-4\iota_0(x),
\end{equation}
and, therefore, $\frv$ is a submodule. This shows too that
$[\iota_0(1),\iota_0(x)]\ne 0$ for any $0\ne x\in\calA$ with $\bar
x=-x$. Moreover,
\[
\begin{split}
[\iota_2(1),\iota_2(x)]&=[[\iota_0(1),\iota_1(1)],\iota_2(x)]\\
  &=[[\iota_0(1),\iota_2(x)],\iota_1(1)]+
    [\iota_0(1),[\iota_1(1),\iota_2(x)]]\\
  &=-[\iota_1(\bar x),\iota_1(1)]+[\iota_0(1),\iota_0(\bar x)]\\
  &=-[\iota_0(1),\iota_0(x)]-[\iota_1(1),\iota_1(x)].
\end{split}
\]
Therefore, $\sum_{i=0}^2[\iota_i(1),\iota_i(x)]=0$, and the
dimension of $\frv$ is at most $5$. But \eqref{eq:iota1iotaxskew},
\eqref{eq:iota1iotaxskewbis} and their analogues for $i=0,1,2$ show
that $[\iota_0(1),\iota_0(x)]$ and $[\iota_1(1),\iota_1(x)]$ are
linearly independent elements of $\frg_{(\bar 0,\bar 0)}$. The
outcome is that the dimension of $\frv$ is $5$. Besides, the
assignment $\iota_i(x)\mapsto G_i$, $[\iota_i(1),\iota_i(x)]\mapsto
-2(H_{i+1}-H_{i+2})$ ($G_i$'s and $H_i$'s as in \eqref{eq:HGDs})
shows that $\frv$ is isomorphic to the irreducible module $\frh$.
Therefore, $\oplus_{i=0}^2\iota_i(\calA)$ is contained in a sum of
irreducible modules for $\frs$ isomorphic either to the adjoint
module or to $\frh$.

Now take any element $0\ne d\in\frg_{(\bar 0,\bar 0)}$ and let
$\calU=\calU(\frs)$ be the universal enveloping algebra of $\frs$.
The $\frs$-module generated by $d$ is $\calU
d=kd+\sum_{i=0}^2\calU[d,\iota_i(1)]$. Let $x_i\in\calA$ be the
elements such that $[d,\iota_i(1)]=\iota_i(x_i)$ ($i=0,1,2$). Then
$\calU d=kd+\sum_{i=0}^2\calU\iota_i(x_i)$. But the sum
$\sum_{i=0}^2\calU\iota_i(x_i)$ is a finite sum of irreducible
modules, each of them isomorphic either to the adjoint module or to
$\frh$, and hence, by complete reducibility, to a finite direct sum
of irreducible modules of these types. Therefore, to prove that
$\calU d$ is a sum of irreducible $\frs$-modules which are either
trivial, adjoint or isomorphic to $\frh$, it is enough to prove that
if $M$ is a module for $\frs$ and $N$ a submodule of $M$ with $N$
either adjoint or isomorphic to $\frh$, and the dimension of $M/N$
is $1$, then $M$ contains a one dimensional submodule (which
necessarily complements $N$). But the Casimir element
$D_0^2+D_1^2+D_2^2\in\calU$ acts as $-2Id$ on the adjoint module,
$-6Id$ on $\frh$ and trivially on the trivial module. Hence the one
dimensional submodule of $M$ sought for is the kernel of the action
of the Casimir element.
\end{proof}

\medskip

By means of \eqref{eq:S4actiononTCJ}, actions of the group $S_4$ on
the exceptional Lie algebras were considered. (Note that the simple
Lie algebra of type $G_2$ appears simply as $\der C$). The previous
Theorem makes easy to embed $S_4$ in the group of automorphisms of
classical Lie algebras.

\begin{examples} Consider the module $W$ for $S_4$ in
\eqref{eq:S4W}.
\begin{romanenumerate}
\item \emph{Orthogonal Lie algebras:}\quad
The module $W$ is endowed with a nondegenerate symmetric
bilinear form $b$ invariant under the action of $S_4$:
$b(w_i,w_j)=\delta_{ij}$ for any $i,j=0,1,2$. Let $(U,b')$ be any
vector space endowed with a nondegenerate symmetric bilinear form.
Then the orthogonal Lie algebra of the orthogonal sum $(W\oplus
U,b\perp b')$ decomposes as:
\[
\frso(W\oplus U,b\perp b')=\frso(W,b)\oplus (W\otimes
U)\oplus\frso(U,b'),
\]
where $\frso(W,b)=\frso_3$ (respectively $\frso(U,b')$) is
identified to the subalgebra of $\frso(W\oplus U,b\perp b')$ which
preserves $W$ (resp.~$U$) and annihilates $U$ (resp.~$W$), and for
any $w\in W$ and $u\in U$, $w\otimes u$ is identified to the linear
map determined by $w'\mapsto b(w,w')u$, $u'\mapsto -b'(u,u')w$, for
any $w'\in W$ and $u'\in U$.

As a module for $\frso_3$, $W$ is isomorphic to the adjoint module,
so $W\otimes U$ is a direct sum of copies of the adjoint module,
while $\frso(U,b')$ is a trivial module. Hence, according to Theorem
\ref{th:g135}, $\frso(W\oplus U,b\perp b')$ is endowed with an
action of $S_4$ by automorphisms.

\item \emph{Special Lie algebras:}\quad Let $U$ be now any vector
space. Then the special linear Lie algebra $\frsl(W\oplus U)$
decomposes as
\[
\frsl(W\oplus U)=\frsl(W)\oplus (W\otimes U^*)\oplus (W^*\otimes
U)\oplus \frgl(U)
\]
with natural identifications. But as in \eqref{eq:glWso3hz},
$\frsl(W)$ decomposes as $\frsl(W)=\frso_3\oplus\frh$ and, as a
module for $\frso_3$, $W$ and $W^*$ are both isomorphic to the
adjoint module. Then $\frsl(W\oplus U)$ decomposes, as a module for
$\frso_3$, as a direct sum of copies of the adjoint module, of
$\frh$ (just one copy) and of the trivial module, so $\frsl(W\oplus
U)$ (or $\frsl_n$ for $n\geq 3$) is endowed with an action of $S_4$
by automorphisms.

\item \emph{Symplectic Lie algebras:}\quad Let now $(U,B')$ be a
vector space endowed with a nondegenerate skew symmetric bilinear
form. Also $W\oplus W^*$ is endowed with the natural skew symmetric
bilinear form $B$, where $W$ and $W^*$ are isotropic subspaces and
$B(f,w)=f(w)$ for any $f\in W^*$ and $w\in W$. The symplectic Lie
algebra of the orthogonal sum $\bigl((W\oplus W^*)\oplus U,B\perp
B'\bigr)$ decomposes as
\[
\begin{split}
\frsp\bigl((W&\oplus W^*)\oplus U,B\perp B'\bigr)\\
&=\frsp(W\oplus W^*,B)\oplus \bigl((W\oplus W^*)\otimes
U\bigr)\oplus \frsp(U,B').
\end{split}
\]
But $\frgl(W)$ is naturally embedded in $\frsp(W\oplus W^*,B)$ as
the subalgebra that leaves both $W$ and $W^*$ invariant. Hence
$\frso_3$, which is contained in $\frgl(W)$, embeds in
$\frsp(W\oplus W^*,B)$ which, as a module for $\frso_3$ is the
direct sum of $\frso_3$, three copies of $\frh$ and three copies of
$\frz$. Again this shows that $\frsp\bigl((W\oplus W^*)\oplus
U,B\perp B'\bigr)$ is endowed with an action of $S_4$ by
automorphisms.
\end{romanenumerate}
\end{examples}

\medskip

\begin{remark}\label{re:Seligman}
The Lie algebras over a field of characteristic $0$ containing a
three dimensional simple Lie algebra $\frs$ such that, as modules
for $\frs$, are a direct sum of copies of the adjoint, the unique
five dimensional irreducible module and the trivial module have been
thoroughly studied in \cite[Chapter 7]{Seligman}.\quad\qed
\end{remark}

\begin{remark}\label{re:TitsTKK}
The Lie algebras $\frg$ whose Lie algebras of derivations contain a
subalgebra isomorphic to $\frso_3$ and such that, as modules for
this subalgebra, they are a direct sum of irreducible modules
isomorphic either to the adjoint module $\frso_3$, the five
dimensional module $\frh$ or the trivial one dimensional module
$\frz$, can be shown to admit a group of automorphisms isomorphic to
$S_4$ exactly as in the proof of Theorem \ref{th:g135}. In
particular, if $J$ is any Jordan algebra (not necessarily unital)
and $\frd$ is a Lie subalgebra of $\der(J)$ containing the inner
derivations, then the Lie algebra $\frg=\bigl(\frso_3\otimes
J\bigr)\oplus\frd$, where $\frd$ is a subalgebra and the bracket is
determined by (see \cite{Tits62}) $[A\otimes x,B\otimes
y]=[A,B]\otimes xy +\frac{1}{2}\trace(AB)[L_x,L_y]$ and
$[d,(A\otimes x)]=A\otimes d(x)$ for any $A,B\in\frso_3$, $x,y\in
J$, and $d\in\frd$, is a Lie algebra satisfying the conditions
above. Note that in case $-1\in k^2$, then $\frso_3$ is isomorphic
to $\frsl_2$, and the construction above becomes the well-known
Tits-Kantor-Koecher construction $TKK(J)$ (see \cite[Example
3.2]{EOJAlgebra}). In particular, if $J$ is the Jordan superalgebra
of type $D(t)$ or $F$, then this construction will give the Lie
superalgebra of type $D(2,1;t)$ or $F(4)$ respectively.\quad\qed
\end{remark}

\medskip

Two more comments are in order here. In a previous paper
\cite{EOJAlgebra}, the authors have shown how to define an action of
$S_4$ on the Lie algebra $\calK(\calA,-)$ attached to a structurable
algebra $(\calA,-)$ by means of Kantor's construction
\cite{Allisonmodelsisotropic} in case $-1$ is a square on the ground
field. The previous Theorem provides a natural interpretation: The
Lie algebra $\calK(\calA,-)$ contains a subalgebra isomorphic to
$\frsl_2$ such that, as a module for $\frsl_2$, $\calK(\calA,-)$ is
a direct sum of copies of $\frsl_2$, of its five dimensional
irreducible module in $\frgl(\frsl_2)$ and the trivial module.
However, if $-1$ is a square, then $\frsl_2$ is isomorphic to
$\frso_3$ and, after identifying $\frsl_2\simeq \frso_3$, the five
dimensional irreducible module is the module $\frh$ considered so
far.

\smallskip

Also, the Lie algebras containing a subalgebra isomorphic to
$\frsl_2$ such that, as a module for $\frsl_2$, are a direct sum of
copies of $\frsl_2$, of its five dimensional irreducible module in
$\frgl(\frsl_2)$ and the trivial module are, essentially, the
$BC_1$-graded Lie algebras of type $B_1$ (see \cite{BenkartSmirnov})
and the references therein). These Lie algebras present a
decomposition as in \eqref{eq:gso3hd}:
\[
\frg=\bigl(\frsl_2\otimes\calH\bigr)\oplus
\bigl(\frh\otimes\calS\bigr)\oplus \frd.
\]
Take the standard basis $\{e,f,h\}$ of $\frsl_2$ with $[e,f]=h$,
$[h,e]=2e$, $[h,f]=-2f$. The action of $\ad h$ gives a $5$-grading:
$\frg=\frg_{-2}\oplus\frg_{-1}\oplus\frg_0\oplus\frg_1\oplus\frg_2$,
where $\frg_i=\{x\in\frg: [h,x]=ix\}$ ($i=-2,-1,0,1,2$). There is
just one extra condition in the definition of the $BC_1$-graded Lie
algebras of type $B_1$:
$\frg_0=[\frg_{-2},\frg_2]+[\frg_{-1},\frg_1]$. With the notations
as in the proof of Theorem \ref{th:g135}, this is equivalent to the
condition $\frd=d_{\calH,\calH}+d_{\calS,\calS}$.

\providecommand{\bysame}{\leavevmode\hbox
to3em{\hrulefill}\thinspace}
\providecommand{\MR}{\relax\ifhmode\unskip\space\fi MR }
\providecommand{\MRhref}[2]{%
  \href{http://www.ams.org/mathscinet-getitem?mr=#1}{#2}
} \providecommand{\href}[2]{#2}


\begin{thebibliography}{ZSSS82}



\bibitem[All78]{All78}
B.~N. Allison, \emph{A class of nonassociative algebras with
involution
  containing the class of {J}ordan algebras}, Math. Ann. \textbf{237} (1978),
  no.~2, 133--156. 

\bibitem[All79]{Allisonmodelsisotropic}
\bysame, \emph{Models of isotropic simple {L}ie algebras}, Comm.
Algebra
  \textbf{7} (1979), no.~17, 1835--1875.

\bibitem[BS03]{BartonSudbery}
C.~H. Barton and A.~Sudbery, \emph{Magic squares and matrix models
of {L}ie
  algebras}, Adv. Math. \textbf{180} (2003), no.~2, 596--647.

\bibitem[BSm03]{BenkartSmirnov}
G.~Benkart and O.~Smirnov,
     \emph{Lie algebras graded by the root system {$\rm BC\sb 1$}},
   J. Lie Theory \textbf{13} (2003), no.~1, 91--132.

\bibitem[BZ96]{BZ96}
G.~Benkart and E.~Zelmanov, \emph{Lie algebras graded by finite root
systems and intersection matrix algebras}, Invent. Math.
\textbf{126}, (1996), no.~1, 1--45.

\bibitem[BE03]{BE03}
G.~Benkart and A.~Elduque, \emph{The Tits construction and the
exceptional simple classical Lie superalgebras}, Q. J. Math.
\textbf{54} (2003), no.~2, 123--137.

\bibitem[Eld04]{Eld04}
A.~Elduque, \emph{The magic square and symmetric compositions}, Rev.
Mat. Iberoamericana \textbf{20} (2004), no.~2, 475--491.

\bibitem[EO00]{EOMathZ}
A.~Elduque and S.~Okubo, \emph{On algebras satisfying {$x\sp 2x\sp
2=N(x)x$}}, Math. Z. \textbf{235} (2000), no. 2, 275--314.

\bibitem[EO06]{EOJAlgebra}
\bysame, \emph{Lie algebras with $S_4$-action and structurable
algebras}, J.~Algebra, to appear.

\bibitem[FH91]{FultonHarris}
W.~Fulton and J.~Harris, \emph{Representation theory. A first
course}, Graduate Texts in Mathematics \textbf{129},
    Springer-Verlag, New York, 1991.

\bibitem[Hum78]{Humphreys}
J.E.~Humphreys, \emph{Introduction to {L}ie algebras and
representation theory}, Graduate Texts in Mathematics \textbf{9},
      Springer-Verlag, New York, 1978.

\bibitem[Jac68]{JacobsonJordan}
Nathan Jacobson, \emph{Structure and representations of {J}ordan
algebras},
  American Mathematical Society Colloquium Publications, Vol. XXXIX, American
  Mathematical Society, Providence, R.I., 1968. 

\bibitem[LM02]{LandsbergManivel}
J.~M. Landsberg and L.~Manivel, \emph{Triality, exceptional {L}ie
algebras and
  {D}eligne dimension formulas}, Adv. Math. \textbf{171} (2002), no.~1, 59--85.

\bibitem[Oku05]{OkuboLAA}
S.~Okubo, \emph{Symmetric triality relations and structurable
algebras}, Linear Algebra Appl. \textbf{396} (2005), 189--222.

\bibitem[Sch95]{Schafer}
Richard~D. Schafer, \emph{An introduction to nonassociative
algebras}, Dover
  Publications Inc., New York, 1995. 

\bibitem[Sel88]{Seligman}
George B.~Seligman, \emph{Constructions of {L}ie algebras and their
modules}, Lecture Notes in Mathematics \textbf{1300},
Springer-Verlag, Berlin, 1988.

\bibitem[Tit62]{Tits62}
J.~Tits, \emph{Une classe d'alg\`ebres de {L}ie en relation avec les
              alg\`ebres de {J}ordan}, Nederl. Akad. Wetensch. Proc. Ser. A \textbf{65}
              = Indag. Math. \textbf{24} (1962), 530--535.

\bibitem[Tit66]{Tits66}
\bysame, \emph{Alg\`ebres alternatives, alg\`ebres de {J}ordan et
alg\`ebres de
  {L}ie exceptionnelles. {I}. {C}onstruction}, Nederl. Akad. Wetensch. Proc.
  Ser. A \textbf{69} = Indag. Math. \textbf{28} (1966), 223--237.

\bibitem[ZSSS82]{ZSSS}
K.A.~Zhevlakov, A.M~Slin'ko, I.P.~Shestakov, A.I.~Shirshov,
   \emph{Rings that are nearly associative},
   Pure and Applied Mathematics \textbf{104}, Academic Press Inc.,
   New York, 1982.


\end{thebibliography}
\end{document}